\documentclass[12pt]{article}
\usepackage{latexsym, amsmath, amssymb, amsfonts, amscd}
\usepackage{graphics}
\usepackage[dvips]{epsfig}

\usepackage[pagebackref]{hyperref}

\newtheorem{theorem}{Theorem}
\newtheorem{proposition}{Proposition}
\newtheorem{lemma}{Lemma}

\newtheorem{definition}{Definition}

\newtheorem{question}{Question}

\newcommand{\oper}[2]{\newcommand{#1}{\mathop{\mathrm{#2}}\nolimits} }
\oper{\DGA}{DGA} \oper{\END}{End} \oper{\fol}{{\it fol}}

\newcommand{\call}{{\cal L}}
\newcommand{\calo}{{\cal O}}

\def\qed{\rule{2.3mm}{2.3mm}}

\def\lcf{\lbrack\! \lbrack}
\def\rcf{\rbrack\! \rbrack}

\def\dbar{\overline\partial}

\def\zbar{\overline z}

\def\oL{\overline L}
\def\oom{\overline\omega}

\newcommand{\CC}{\mathbb{C}}

\def\cp{{\bf C}{\bf P}}

\def\cl{{\cal L}}

\def\co{{\cal O}}

\def\dbar{\overline\partial}

\def\zbar{\overline z}

\def\osigma{\overline\sigma}

\def\oomega{\overline\omega}

\def\osigma{\overline\sigma}

\def\oz{\overline{z}}

\newcommand{\lie}[1]{\mathfrak{#1}}

\newcommand{\lbra}[2]{\lcf #1, #2 \rcf}

\newcommand{\pair}[2]{\langle #1, #2 \rangle}

\newcommand{\bproof}{\noindent{\it Proof: }}
\newcommand{\eproof}{\hfill \qed \vspace{0.2in}}

\begin{document}
\title{\bf Differential Gerstenhaber Algebras \\
of  Generalized Complex Structures}
\author{
Daniele Grandini,
\thanks{ Address:
Department of Mathematics, University of California at Riverside,
    Riverside CA 92521, U.S.A.} \
Yat-Sun  Poon, \thanks{ Address:
    Department of Mathematics, University of California at Riverside,
    Riverside CA 92521, U.S.A., Email: ypoon@ucr.edu. Partially
    supported by NSF DMS-0906264 and Mathematical Sciences Center of Tsinghua
    University} \
   Brian Rolle \thanks{Address:
 Department of Mathematics, University of California at Riverside,
    Riverside CA 92521, U.S.A.}
  }
\date{ August 26, 2011}
\maketitle
\begin{abstract} Associated to every generalized complex structure is a differential Gerstenhaber algebra ($\DGA$).
When the generalized complex structure deforms, so does the associated $\DGA$.
In this paper, we identify the infinitesimal conditions when the $\DGA$ is invariant as the generalized complex
structure deforms. We prove that the infinitesimal condition is always integrable. When the underlying manifold is a holomorphic
Poisson nilmanifolds, or simply a group in the general, and the geometry is invariant, we find a general construction
to solve the infinitesimal conditions under some geometric conditions.
Examples and counterexamples of existence of solutions to the infinitesimal
conditions are given.
\end{abstract}


\section{Introduction}

A few years ago,  the second author computed the weak Frobenius structure on the moduli space of
the Barannikov-Kontsevich's extended deformation \cite{BK} of
the complex structure on a primary Kodaira surface \cite{Poon}.
Among other observations, one could see from \cite[Table (45)]{Poon} that the restriction of the weak
Frobenius structure to the even part of the extended moduli space is trivial. The parameter space of the
even part of the extended moduli at the unperturbed point is contained in
\begin{equation}
\oplus_{k={\mbox{even}}}H^k_J, \quad \mbox{ where } \quad H^k_J=\oplus_{p+q=k}H^q(M,\wedge^p T^{1,0})
\end{equation}
and $T^{1,0}$ is the holomorphic tangent bundle  of the complex manifold $M$.
The computation in \cite{Poon} dwells in the fact that the primary Kodaira surface was chosen to be
a nilmanifold and the
complex structure he worked with is invariant. Along the line of thoughts in \cite{Nomizu} \cite{GMPP}, the
Dolbeault cohomology could be computed by means of algebraic methods.

Thanks to the work of Hitchin \cite{Hitchin} and Gualtieri \cite{Marco}, it is now well known that the degree-2
portion of the extended deformation is realized by deformation of \it generalized geometry. \rm While we will
provide further details on generalized geometry in Section \ref{sec:gcx}, at this stage we simply note
that  the parameter
space of generalized deformation is the degree-2 portion of extended deformation.
\begin{equation}\label{degree 2}
H_J^2=H^0(M, \wedge^2T^{1,0})\oplus H^1(M, T^{1,0})\oplus H^2(M, \calo)
\end{equation}
where $\calo$ is the structure sheaf of the complex manifold $M$.

The key ingredient in constructing the weak Frobenius structure on extended deformation is a variation of the
exterior product structure when the concerned cohomology spaces vary. However, it is also known that the
differential geometric object controlling the extended deformations is the differential Gerstenhaber algebras
$(\DGA)$
associated to each (extended) complex structure \cite{CZ} \cite{Mer} \cite{Mer-note} \cite{Zhou}.  We will provide necessary details on the construction of $\DGA$s in Section \ref{sec:DGA}.
  This structure contains the exterior differential algebra as a sub-structure.
In this context, we could paraphrase a result of \cite{Poon} in a context of generalized complex geometry, and
say that the exterior differential algebras along a generalized deformation of a
 primary Kodaira surface is rigid, meaning that all the exterior differential algebras are quasi-isomorphic to the unperturbed one. From this perspective, we seek a general understanding of the rigidity of the full differential Gerstenhaber algebra structures.

\begin{question}\label{general question} Suppose that $M$ is a manifold with generalized complex structure $J$. Let $\DGA(0)$ be the
associated differential Gerstenhaber algebra. Suppose that $\Gamma(t)$ is a family of deformation of $J$
along  generalized complex structure with parameter $t$,
with associated differential Gerstenhaber algebra $\DGA(t)$. Under what condition
will $\DGA(t)$ be quasi-isomorphic to $\DGA(0)$?
\end{question}

The infinitesimal counter-part of $\Gamma(t)$ is $\Gamma_1$, which represents an element in the cohomology
space $H^2_J$. If there is a quasi-isomorphism $\Phi(t)$, depending on $t$, we consider its infinitesimal version $\phi$. The pair $\Gamma_1$ and $\phi$ will be addressed as \it compatible pair. \rm Together, they have to satisfy
a set of constraints as given in Definition \ref{compatible pair}. The main result in Section \ref{sec:DGA} is Theorem
\ref{thm:integrable}, which states essentially that compatible pairs are always integrable. Therefore, answers to
Question \ref{general question} above are reduced to infinitesimal level.

In identity (\ref{degree 2}), we see that there are three special kinds of deformations to analyze. Those from
$H^1(M, T^{1,0})$ are due to classical complex deformation theory. Those from
$H^2(M, \calo)$ are due to $B$-field transformations if the underlying complex structure is K\"ahlerian
\cite{Marco}. Therefore, we focus on those in the component $H^0(M, \wedge^2T^{1,0})$. As we will explain later,
this class of deformation is due to holomorphic Poisson structures, objects under investigation
from various perspectives \cite{Goto} \cite{Hitchin3} \cite{Hitchin2}. If the holomorphic Poisson structure
has full rank everywhere, it leads to a deformation from a classical complex structure $J$
to a symplectic structure $\Omega$. If the induced differential Gerstenhaber algebras along this deformation is
rigid, then $\DGA(J)$ and $\DGA(\Omega)$ are quasi-isomorphic. It presents the complex manifold
$(M,J)$ and the symplectic manifold $(M, \Omega)$ as a weak mirror pair in the sense of Merkulov
\cite{Mer-note}. An investigation on such possibility also motivates this paper. Therefore, in Section
\ref{sec:Poisson} we refine our analysis in Section \ref{sec:DGA} to holomorphic Poisson manifolds, and illustrate our theory with a computation on a Hopf surface.

For nilmanifolds, i.e. the compact quotient of simply connected nilpotent Lie groups, it is known for
a very long time that
the DeRham cohomology is given by invariant elements \cite{Nomizu}. From our current perspective,
the invarant $\DGA$ with an invariant symplectic structure on a nilmanifold is quasi-isomorphic to the
full $\DGA$ of the symplectic structure.
For a large class of nilmanifolds examples, we also know that the invariant $\DGA$ theory for invariant complex
structures is  quasi-isomorphic to the $\DGA$ of the corresponding nilmanifolds \cite{CFGAU} \cite{GMPP}
  \cite{Poon} \cite{Rolle}. Therefore,  we reduce the theory
  in the previous sections in terms of invariant objects on Lie algebras and develop a method to
  construct compatible pairs on a class of holomorphic Poisson algebras in Section \ref{sec:algebra}.
  Finally, in Section \ref{sec:example} we analyze all non-trivial real four-dimensional examples.
  Among other observations, we conclude that the differential Gerstenhaber algebra structures are rigid
  when one deforms the complex structure on a Kodaira surface by a holomorphic Poisson structure.
  It extends the results in \cite{Poon} on weak Frobenius structures, at least along the degree-2 direction
  of the extended moduli space. On the other hand, we also discover an example of holomorphic symplectic
  algebra on which there is no compatible pair. Therefore, a solution to Question \ref{general question} is
  non-trivial.

In this notes, we assume that readers are familiar with the concepts of
Lie algebroids and Lie bialgebroids. Otherwise,
\cite{LWX} and \cite{Mac-2} are our references. On Differential Gerstenhaber
algebras, we rely on \cite{Mac-2} and \cite{Poon} for their formal aspects.
For generalized complex structures, our references are \cite{Hitchin} and \cite{Marco}. Much of the computation
in Section \ref{sec:algebra} and Section \ref{sec:example} could be found in the third author's thesis.
Therefore, our presentation will be relatively sketchy.

 \section{Generalized complex structures}\label{sec:gcx}

Let $M$ be a smooth connected manifold without boundary. Denote its tangent and cotangent bundle
respectively by $T$ and $T^*$. If $V$ is a vector bundle on $M$, we denote its space of sections
by $C^\infty(V)$.  Generic vector fields will be denoted by $X$ and $Y$. One-forms are denoted
by $\alpha$ and $\beta$.  On the bundle  $T\oplus T^*$, there is a natural pairing defined by
\begin{equation}\label{pairing}
\pair{X+\alpha}{Y+\beta}=\frac12(\alpha(Y)+\beta(X)).
\end{equation}
As this pairing is non-degenerate, it identifies the bundle $T\oplus T^*$ to its dual.
We choose the identification to be
\begin{equation}\label{identification}
\sigma : T\oplus T^*\to (T\oplus T^*)^*,
\quad
\sigma(X+\alpha)(Y+\beta)=2\pair{X+\alpha}{Y+\beta}.
\end{equation}
The Courant bracket \cite{LWX} is the real bilinear map on $C^\infty(T\oplus T^*)$ defined by
\begin{equation}\label{courant}
\lbra{X+\alpha}{Y+\beta}=
[X, Y]+{\call}_X\beta- {\call}_Y\alpha-
 \frac{1}{2}d(\iota_X\beta- \iota_Y\alpha).
 \end{equation}
 The Courant bracket, the non-degenerate pairing above, together with the natural projection on the tangent component
 make $T\oplus T^*$ a standard example of a Courant algebroid
 \cite{Courant} \cite{LWX}.

 An almost generalized complex structure is a real bundle map
 $J: T\oplus T^*\to T\oplus T^*$ such that $J\circ J=-$identify and $J^*=-J$.
 Let $L$ be the bundle of $+i$-eigenvectors with respect to $J$ and over the
 complex numbers. With respect to the non-degenerate pairing, $L$ is maximal isotropic.
  So is its conjugate bundle $\overline{L}$.
 The choices of the tensorial object $J$ with the given prescription is equivalent to
 the choice of maximal isotropic subbundle $L$ such that
 $L\cap {\overline{L}}$ is trivial \cite{Marco}.

 An almost generalized complex structure is said to be integrable if and only if the space
  $C^\infty(L)$ is closed under the Courant bracket.
 By complex conjugation, it is of course equivalent to $C^\infty(\overline{L})$ being closed.
 In such case, the structure $J$, or equivalently, either the bundle $L$ or the bundle
  $\overline{L}$ is said to be a
 generalized complex structure. It is now well known that complex structures in the classical sense
 are generalized complex. So are symplectic structures. For classical complex structure, the complexified
 tangent bundle splits into the direct sum of type $(1,0)$ and type $(0,1)$ vectors. Their related
 bundles are denoted by $T^{1,0}$ and $T^{0,1}$ respectively. Their dual bundles are denoted by
 $T^{*(1,0)}$ and $T^{*(0,1)}$. Then the corresponding bundles $L$ and $L^*$ are
 \[
 L=T^{1,0}\oplus T^{*(0,1)}, \quad L^*\cong {\overline{L}}=T^{0,1}\oplus T^{*(1,0)}.
 \]

 If $\omega$ is a symplectic form on the manifold $M$, then
 \[
 L=\{ X-i\iota_X\omega: X\in C^\infty(T)\}, \quad L^*\cong {\overline{L}}=\{ X+i\iota_X\omega: X\in C^\infty(T)\}
 \]
 represent an example of a generalized complex structure.

Since $L$ is isotropic, the restriction of the Courant bracket on $L$
makes it a Lie algebroid whenever the generalized complex structure is integrable.
 As such, it has a Lie algebroid differential acting on the exterior algebra of the dual bundle \cite{Mac-2}.
 \begin{equation}
 \partial: C^\infty(\wedge^n L^*)\to C^\infty(\wedge^{n+1}L^*).
 \end{equation}
 Using the identification as given in (\ref{identification}), we identify $\overline{L}=L^*$. Then
 \begin{equation}
 \partial: C^\infty(\wedge^n \overline{L})\to C^\infty(\wedge^{n+1} \overline{L}).
 \end{equation}
Similarly, $\overline{L}\cong L^*$ is also a Lie algebroid. Its Lie algebroid differential
is precisely the conjugation of the above operator:
\begin{equation}
\dbar: C^\infty(\wedge^n L)\to C^\infty(\wedge^{n+1}L).
\end{equation}
As noted in \cite[Theorem 2.6]{LWX}, $(L, {\overline{L}})$ forms a Lie bialgebroid. It means
that for any sections $\ell_1$ and $\ell_2$ of
the bundle $L$,
\begin{equation}\label{compatible}
\dbar\lbra{\ell_1}{\ell_2}=\lbra{\dbar\ell_1}{\ell_2}-\lbra{\ell_1}{\dbar\ell_2}.
\end{equation}

Making use of \cite[Theorem 7.5.2]{Mac-2}, we deduce that the space of sections of the exterior
algebra generated by $L$, $C^\infty(\wedge^\bullet L)$ carries the structure of a differential
Gerstenhaber algebra structures, with the Courant bracket, exterior product and Lie algebroid
differential of $\overline{L}\cong L^*$.
We denote it by
\begin{equation}
DGA(J):= (C^\infty(\wedge^\bullet L), \lbra{-}{-}, \wedge, \dbar).
\end{equation}
In this context, the bracket on $C^\infty(\wedge^\bullet L)$ is known as Schouten bracket \cite{Mac-2}.

The integrability implies that the restriction of the Courant bracket on $C^\infty(L^*)$ satisfies
the Jacobi identity. In terms of the operator
$\dbar$, it is equivalent to $\dbar\circ\dbar=0$. Therefore,
$\dbar: C^\infty(\wedge^n L)\to C^\infty(\wedge^{n+1}L)$ determines a differential
complex, and hence generates cohomology spaces.  i.e. for all $k\geq 1$,
\[
H^k_J=\frac{\ker\dbar:\wedge^kL \to \wedge^{k+1}L}{{\mbox{\rm Image } }\dbar:\wedge^{k-1}L\to\wedge^kL}.
\]
Given the identity (\ref{compatible}),
the cohomology spaces inherit a Gerstenhaber algebra structure.

When the generalized complex structure is classical,  one could verify that if $\oomega$ is
a type $(0,k)$-form, then $\dbar\oomega$ is the classical $\dbar$-operator in complex analysis
on $\CC^n$ \cite{Rolle}. On the other hand, if $Z$ is a $(1,0)$-vector field and $\overline X$ is
a $(0,1)$-vector field, then
\begin{equation}\label{CR}
\dbar_{\overline X}Z=[Z, {\overline X}]^{1,0}.
\end{equation}
This is precisely the Cauchy-Riemann operator \cite{Gau} \cite{Brian}.
The cohomology of degree-$k$ in this case is
\begin{equation}
H^k_J=\oplus_{p+q=k}H^q(M, \wedge^pT^{1,0}).
\end{equation}
Using Dolbeault theory, the elements in these cohomology spaces are represented by $\dbar$-closed
 $(0,q)$-forms with coefficients in holomorphic $(p,0)$-vector fields.

On the other hand, if a generalized complex structure is defined by a symplectic form $\omega$, then
\begin{equation}
\dbar({X}-i\iota_{X}\theta)=-2id\iota_{X}\theta
\end{equation}
for all $X$ in $C^\infty(T_\CC)$ \cite{Brian}. In particular, the k-th cohomology of this complex is the k-th
complexified deRham cohomology of the manifold $M$.

As a subbundle of $(T\oplus T^*)_\CC$, the bundle $L$ has a natural projection  $\rho$
onto the direct summand $T_\CC$.
The type of a generalized complex structure at a point of the manifold $M$ is defined to be
the complex co-dimension
of the projection of $L$ in $T_\CC$ over the concerned point \cite{Marco}. From the description above, one sees that
the type of a classical complex structure on a real $2n$-dimensional manifold is
equal to $n$. All symplectic structures are type-0 generalized complex structures.

\section{Deformation of generalized complex structures}\label{sec:DGA}

A deformation of a generalized complex structure is given by a section $\Gamma$
of $\wedge^2L$ \cite{LWX} \cite{Marco}.
To be more precise,
\begin{equation}
L_{\overline\Gamma}=\{\ell+{\overline\Gamma}(\ell): \ell\in C^\infty(L)\},
\quad \mbox{ and } \quad {\overline L}_{\Gamma}
=\{{\overline\ell}+\Gamma({\overline\ell}): {\overline\ell}\in C^\infty({\overline{L}})\}.
\end{equation}
$L_{\overline\Gamma}\cap {\overline L}_{\Gamma}=\{0\}$ if and only
if $\Gamma\circ {\overline \Gamma}$ does not have non-trivial fixed points \cite{Brian}.
The deformed generalized complex structure $(L_{\overline\Gamma}, {\overline L}_{\Gamma})$
is integrable if and only if
$\Gamma$ satisfies the Maurer-Cartan equation \cite[Theorem 6.1]{LWX}:
\begin{equation}\label{MC equation}
\dbar\Gamma+\frac12 \lbra{\Gamma}{\Gamma}=0.
\end{equation}

The infinitesimal version of the Maurer-Cartan equation is simply $\dbar \Gamma_1=0$.
Therefore, it represents an element in the second cohomology
$H^2_J$ of the differential Gerstenhaber algebra of the unperturbed generalized complex structure $J$.

\subsection{Deformation of associated $\DGA$}
Let $\overline{\delta}$ be the Lie algebroid differential of ${\overline{L}}_{\Gamma}$.
Due to our natural pairing (\ref{pairing}), it acts on
the conjugate bundle $L_{\overline \Gamma}$. Therefore, we have the new differential Gerstenhaber algebra
\begin{equation}
DGA(J_{\Gamma})=(\wedge^\bullet L_{\overline\Gamma}, \lbra{-}{-}, \wedge, \overline{\delta}).
\end{equation}

Meanwhile, for $\Gamma$ sufficiently close to zero, $L$ and  ${\overline{L}}_{\Gamma}$
are also transversal in $(T\oplus T^*)_\CC$.
By \cite[Theorem 2.6]{LWX}, $L$ and ${\overline{L}}_{\Gamma}$ form a Lie bialgebroid.
We could denote the Lie algebroid differential of
the Lie algebroid ${\overline{L}}_{\Gamma}$ acting on $L$ by $\dbar_\Gamma$.
Since $L_{\overline\Gamma}$ is simply the graph of the map $\overline\Gamma$, there is a natural
map from $L$ to $L_\Gamma$. It enables one to identify the differential $\dbar_\Gamma$. The
computation below is a consequence of \cite[Theorem 2.6]{LWX} and \cite[Theorem 6.1]{LWX}.
It should be well known to experts. We outline a proof here
for completeness. A complete proof for a case most relevant to this paper could be found in
\cite{Brian}.

\begin{proposition}\label{reference}
The pair
$L$ and $\overline{L}_{\Gamma}$ forms a Lie bialgebroid. The Lie algebroid
differential $\overline{\partial}_{\Gamma}$ for the
deformed Lie algebroid $\overline{L}_{\Gamma}$ acting on $L$ is
given by $\overline {\partial}+\lbra{\Gamma}{-}.$ i.e. for any
section $\ell$ of $L$,
\[
\overline{\partial}_{\Gamma} \ell  =\overline{\partial}%
\ell+\lbra{\Gamma}{\ell}.
\]
\end{proposition}\label{Daniele 1}
\bproof By definition, the vector bundle $\overline{L}_{\Gamma}$ is maximally isotropic
in the Courant algebroid $(  T\oplus T^*)  _{\mathbb{C}}$. Since $\Gamma$
satisfies the Maurer-Cartan equation, it follows that the space of sections of $\overline{L}_{\Gamma}$
is closed with respect to the Courant bracket. To find a more precise description, we follow
the computation in \cite{LWX}.

For ${\overline\sigma}\in C^\infty({\overline L})$, $\ell_1, \ell_2\in C^\infty(L)$, define a
 Lie derivative by
\[
(\call_{\ell_1}{\overline\sigma})\ell_2=\rho(\ell_1)({\overline\sigma}(\ell_2))-{\overline\sigma}(\lbra{\ell_1}{\ell_2}),
\]
where $\rho$ is the natural projection from $(T\oplus T^*)_\CC$ onto $T_\CC$.
The property of the Lie derivative in algebroid theory could be found in \cite{Mac-2}. Follow
 \cite[Identity (23)]{LWX}, for any
 $\overline{\ell}_1,\overline{\ell}_2\in C^\infty(\overline{L})$ define
 \begin{equation}\label{Daniele 2}\lbra{\overline{\ell}_1}{\overline{\ell}_2}_{\Gamma}
={\mathcal{L}}_{\Gamma\overline{\ell}_1}\overline{\ell}_2-{\mathcal{L}}_{\Gamma\overline{\ell}_2}\overline{\ell}_2
+\overline{\partial}(\overline{\ell}_1(\Gamma\overline{\ell}_2)).
\end{equation}
As noted in the proof of \cite[Theorem 6.1]{LWX},  $\Gamma$ satisfies the
Maurer-Cartan equation if and only if
\[
\lbra{\overline{\ell}_1+\Gamma\overline{\ell}_1}{\overline{\ell}_2+\Gamma\overline{\ell}_2}
=\lbra{\overline{\ell}_1}{\overline{\ell}_2}+\lbra{\overline{\ell}_1}{\overline{\ell}_2}_{\Gamma}
+\Gamma\left(\lbra{\overline{\ell}_1}{\overline{\ell}_2}+\lbra{\overline{\ell}_1}{\overline{\ell}_2}_{\Gamma}\right).
\]

Now we are ready to compute the Lie algebroid differential of $\overline{L}_{\Gamma}$ with $L$ as its dual.
For every $\ell\in C^\infty(L)$ and $\ell_1,\ell_2\in C^\infty(\overline{L})$,
\begin{eqnarray*}
&&\left(\dbar_\Gamma\ell\right)(\overline{\ell}_1+\Gamma\overline{\ell}_1,\overline{\ell}_2+\Gamma\overline{\ell}_2)\\
&=&\rho(\overline{\ell}_1+\Gamma\overline{\ell}_1)(\ell(\overline{\ell}_2))
-\rho(\overline{\ell}_2+\Gamma\overline{\ell}_2)(\ell(\overline{\ell}_1))
-\ell(\lbra{\overline{\ell}_1+\Gamma\overline{\ell}_1}{\overline{\ell}_2+\Gamma\overline{\ell}_2})\\
&=&\rho(\overline{\ell}_1+\Gamma\overline{\ell}_1)(\ell(\overline{\ell}_2))
-\rho(\overline{\ell}_2+\Gamma\overline{\ell}_2)(\ell(\overline{\ell}_1))\\
&&-\ell\left(\lbra{\overline{\ell}_1}{\overline{\ell}_2}+\lbra{\overline{\ell}_1}{\overline{\ell}_2}_{\Gamma}
+\Gamma\left(\lbra{\overline{\ell}_1}{\overline{\ell}_2}+\lbra{\overline{\ell}_1}{\overline{\ell}_2}_{\Gamma}\right)\right).
\end{eqnarray*}
Since the image of a section in $\overline L$ under $\Gamma$ is a section
of $L$, and $L$ isotropic, the above is equal to
\begin{eqnarray*}
&=&\rho(\overline{\ell}_1+\Gamma\overline{\ell}_1)(\ell(\overline{\ell}_2))-\rho(\overline{\ell}_2
+\Gamma\overline{\ell}_2)(\ell(\overline{\ell}_1))
-\ell\left(\left(\lbra{\overline{\ell}_1}{\overline{\ell}_2}+\lbra{\overline{\ell}_1}{\overline{\ell}_2}_{\Gamma}\right)\right)\nonumber\\
&=&\left(\overline{\partial}{\ell}\right)(\overline{\ell}_1,\overline{\ell}_2)+\rho(\Gamma\overline{\ell}_1)(\ell(\overline{\ell}_2))
-\rho(\Gamma\overline{\ell}_2)(\ell(\overline{\ell}_1))-\ell\left(\lbra{\overline{\ell}_1}{\overline{\ell}_2}_{\Gamma}\right).
\end{eqnarray*}
The proof of this proposition is completed if we could show that
\[
\lbra{\Gamma}{\ell}(\overline{\ell}_1,\overline{\ell}_2)
=\rho(\Gamma\overline{\ell}_1)(\ell(\overline{\ell}_2))
-\rho(\Gamma\overline{\ell}_2)(\ell(\overline{\ell}_1))
-\ell\left(\lbra{\overline{\ell}_1}{\overline{\ell}_2}_{\Gamma}\right).
\]
It is now a matter of definition of Lie derivative to show that the right hand side of the above is equal to
\[
\overline{\ell}_2\left(\lbra{\Gamma\overline{\ell}_1}{\ell}\right)
-\overline{\ell}_1\left(\lbra{\Gamma\overline{\ell}_2}{\ell}\right)
-\rho(\ell)\left(\overline{\ell}_1(\Gamma\overline{\ell}_2)\right).
\]
Finally, the following identity always hold for any section $\Gamma$
of $\wedge^2 L$, $\ell$ of $L$ and ${\overline\ell}_1
,{\overline\ell}_2$ of $\overline L$.
\begin{equation}\label{21 identity}
\lbra{\Gamma}{\ell}(\overline{\ell}_1,\overline{\ell}_2)
=\overline{\ell}_2\left(\lbra{\Gamma\overline{\ell}_1}{\ell}\right)
-\overline{\ell}_1\left(\lbra{\Gamma\overline{\ell}_2}{\ell}\right)
-\rho(\ell)\left(\overline{\ell}_1(\Gamma\overline{\ell}_2)\right)
\end{equation}
See \cite{Brian} for a detailed proof for (\ref{21 identity}).
Therefore, the proof of the proposition is completed.
\eproof

 As far as analyzing the deformation of
  associated differential Gerstenhaber algebras is concerned,
  the above observation reduces an analysis to one on deformation of differentials.
  As the bundle structure could remain constant, we now focus on the variation
  from $\dbar$ to $\dbar_\Gamma$, and hence denote the respective
  differential Gerstenhaber algebras by $\DGA(\dbar)$ and $\DGA(\dbar_\Gamma)$.

To compare $\DGA(\dbar_\Gamma)$ with $\DGA(\dbar)$, we identify the constraints for them to be homomorphic.
If the homomorphism induces an isomorphism at cohomology level, these two $\DGA$s are said to be quasi-isomorphic.
In such case, the generalized complex structure $J$ and the deformed one are also said to form a weak mirror pair
\cite{Mer-note}, \cite{CLP}, \cite{COP}. Let
\[
\Phi: L \to L
\]
be a vector bundle homomorphism depending on $\Gamma$.
It induces a homomorphism of the exterior
algebra generated by $L$. That is
\begin{equation}\label{asso homo}
\Phi(A\wedge B):=\Phi(A)\wedge \Phi(B)
\end{equation}
for any $A, B\in C^\infty(\wedge^\bullet L)$.

It is a homomorphism of graded Lie algebras if
for any $A, B$ in $C^\infty(\wedge^\bullet L)$,
\begin{equation}\label{lie homo}
\lbra{\Phi(A)}{\Phi(B)}=\Phi(\lbra{A}{B}).
\end{equation}
If in addition, when the following diagram is commutative for all $0\leq k\leq n$,
\begin{equation}\label{gauge}
\begin{array}
[c]{ccc}
\wedge^kL & \overset{\overline{\partial}_{\Gamma}}{\longrightarrow} & \wedge^{k+1}L\\
\Phi\downarrow &  & \downarrow \Phi\\
\wedge^kL & \overset{\dbar}{\longrightarrow} & \wedge^{k+1}L,
\end{array}
\end{equation}
we have  a
homomorphism of differential Gerstenhaber algebras
\[
\Phi: (\wedge^\bullet L, \lbra{-}{-}, \wedge, \dbar_\Gamma) \to
(\wedge^\bullet L, \lbra{-}{-}, \wedge, \dbar).
\]

In general, if the bundle homomorphism intertwines the differentials as in
(\ref{gauge}) and satisfies (\ref{asso homo}), then the map $\Phi$ is an exterior
differential algebra homomorphism.
If the bundle homomorphism intertwines the differentials as in
(\ref{gauge}) and satisfies (\ref{lie homo}), the map $\Phi$ is
a graded differential algebra homomorphism.

The diagram (\ref{gauge}) above is equivalent to
\begin{equation}
\Phi\circ \dbar_\Gamma= \dbar \circ \Phi.
\end{equation}
By Proposition \ref{Daniele 1}, that means for any section $A$
of the bundle $\wedge^\bullet L$,
\begin{equation}\label{intertwine}
\Phi(\dbar A+\lbra{\Gamma}{A})=\dbar(\Phi A).
\end{equation}

Suppose that $\Omega$ is a closed 2-form and $L$ is a generalized complex structure, then
\[
{\overline L}_{\Omega}=\{X+\alpha+\iota_X\Omega: X+\alpha\in {\overline L}\}
\]
is again a generalized complex structure. This is because the closedness of $\Omega$
and its skew-symmetry imply that
the map
\[
X+\alpha \mapsto X+\alpha+\iota_X\Omega
\]
is an automorphism of the Courant bracket $\lbra{-}{-}$ in (\ref{courant}) and the
non-degenerate bilinear pairing
$\langle -, -\rangle$ in (\ref{pairing}). This map is known as a $B$-field transformation
by the closed 2-form $\Omega$.
It naturally defines an isomorphism between the $\DGA$ of the bundle $L$ and the $\DGA$ of
the $B$-field transformation of $L$.

On infinitesimal level,  deformation theory of generalized complex structures
is elliptic and a copy of Kuranishi theory follows \cite{Marco}. In this spirit of Kuranishi theory,
let $\phi$ be the infinitesimal version of $\Phi$. It is now an endomorphism from the vector bundle
$L$ to $L$. The infinitesimal version of (\ref{intertwine}) becomes
\begin{equation}
\phi(\dbar A)+\lbra{\Gamma_1}{A}=\dbar(\phi A),
\end{equation}
where $\Gamma_1$ is the infinitesimal deformation, representing an element in the
second cohomology $H^2_J$.
The infinitesimal version of (\ref{lie homo}) is
\begin{equation}
\lbra{\phi  A}{B}+\lbra{A}{\phi B}=\phi\lbra{A}{B}.
\end{equation}
In these expressions, the map $\phi$ is an infinitesimal version of
a homomorphism of exterior algebras. Therefore, it is an
endomorphism with following property.
\begin{equation}
\phi(A\wedge B)=(\phi A)\wedge B+A\wedge (\phi B).
\end{equation}

Now we treat $\phi$ as an element in $L^*\otimes L=\END (L)$. More
generally, it is a section of $\END (\wedge^n L)$. On the other hand,
$\Gamma_1$ is a section of $\wedge^2L$.

Now we summarize the above discussion with a concept and its implication
to variation of the structure of associated differential Gerstenhaber algebras.

\begin{definition}\label{compatible pair} Suppose that $M$ is a manifold with a generalized
complex structure $J$, whose $+i$-eigenbundle is $L$.
A section
$\Gamma_1\in C^\infty(\wedge^2L)$ and a section $\phi\in C^\infty(L^*\otimes L)$ form
a compatible pair if $\dbar\Gamma_1=0$ and
\begin{eqnarray}
 &\dbar(\phi
A)-\phi(\dbar
A)=\lbra{\Gamma_1}{A};&\label{infinitesimal 1}\\
&\lbra{\phi  A}{B}+\lbra{A}{\phi B}=\phi\lbra{A}{B};&\label{phi
endo}\\
&\phi(A\wedge B)=(\phi A)\wedge B+A\wedge (\phi B).&\label{wedge endo}
\end{eqnarray}
\end{definition}

Note that if $\Gamma_1$ is in the center of the Gerstenhaber algebra, i.e.
$\lbra{\Gamma_1}{A}=0$ for all $A\in C^\infty(\wedge^\bullet L)$, then $\phi=0$ is
an obvious solution for the above three identities. For future reference, we note the following

\begin{proposition}\label{central} Suppose that $M$ is a manifold with a generalized
complex structure $J$, whose $+i$-eigenbundle is $L$.
Let
$\Gamma_1\in C^\infty(\wedge^2L)$ be a closed section: $\dbar\Gamma_1=0$.
Then $(\Gamma_1, \phi=0)$ form a compatible pair if and only if $\Gamma_1$ is
central: $\lbra{\Gamma_1}{A}=0$ for all $A\in C^\infty(\wedge^\bullet L)$.
\end{proposition}

\begin{theorem}\label{thm:infinitesimal homo}
Suppose that $M$ is a manifold with a generalized
complex structure $J$, whose $+i$-eigenbundle is $L$. Suppose that
 $\Gamma$ is an integrable deformation with infinitesimal
deformation $\Gamma_1$. If there is a homomorphism $\Phi$ of the bundle
$L$ to itself such that it generates a homomorphism from the
differential Gerstenhaber algebra of the deformation generalized
complex structure to the un-perturbed one, then there exists a
compatible pair $\Gamma_1$ and $\phi$ such that $\dbar\Gamma_1=0$, and up to
first order, $\Gamma$ is equal to $\Gamma_1$ and
$\Phi$ is equal to $1+\phi$.
\end{theorem}

\subsection{Integrability of compatible pairs}

Given a compatible pair $\Gamma_1$ and $\phi$, an immediate issue is whether they
actually come from a deformation $\Gamma$ and a homomorphism of $\DGA$s. In
this section, we apply the principles of
Kuranishi's recursive method to  prove that
 this is the case. We will divide the proof of the following
theorem in several steps.

\begin{theorem}\label{thm:integrable}
Suppose that $M$ is a manifold with a generalized
complex structure $J$, whose $+i$-eigenbundle is $L$. Let $\Gamma_1\in C^\infty(\wedge^2 L)$
and $\phi\in C^\infty(L^*\otimes L)$
be a
compatible pair. Let $t$ be a real variable. Define
\begin{equation}
\Gamma(t)=\sum_{n=1}^\infty (-1)^{n-1}\frac{1}{n!}t^n\phi^{n-1}\Gamma_1, \quad
\Phi(t)=\sum_{n=1}^\infty\frac{1}{n!}t^n\phi^n.
\end{equation}
Then $\Gamma(t)$ satisfies the Maurer-Cartan equation. Moveover, if $\DGA(\Gamma(t))$
represents the Differential Gerstenhaber algebra
$(\wedge^\bullet L, \lbra{-}{-}, \wedge, \dbar_{\Gamma(t)})$, then $\Phi(t)$ is a
homomorphism from $\DGA(\Gamma(t))$ to $\DGA(\Gamma(0))$.
\end{theorem}

The Maurer-Cartan equation at degree-1 with respect to the variable $t$ is
simply $\dbar\Gamma_1=0$. For $n\geq 2$, it is
\begin{equation}
(-1)^{n-1}\frac{1}{n!}\dbar(\phi^{n-1}\Gamma_1)+\frac12\sum_{j+k=n}\lbra{(-1)^{k-1}\frac{1}{k!}\phi^{k-1}\Gamma_1}
{(-1)^{j-1}\frac{1}{j!}\phi^{j-1}\Gamma_1}=0.
\end{equation}
Let the binomial coefficients be
\[
C^n_k=\frac{n!}{k!(n-k)!}.
\]
Then the above equation is equivalent to
\begin{equation}\label{MC n}
\dbar(\phi^{n-1}\Gamma_1)=\frac12\sum_{k=1}^{n-1}C^n_k\lbra{\phi^{k-1}\Gamma_1}{\phi^{n-k-1}\Gamma_1}.
\end{equation}
Similarly, $\Phi(t)$ is a homomorphism of $\lbra{-}{-}$ and $\wedge$ if and
only if for degree $n$, and for any sections
$A$ and $B$ of $ L$,
\begin{equation}
\frac{1}{n!}\phi^n\lbra{A}{B}=\sum_{k+j=n}\lbra{\frac{1}{k!}\phi^kA}{\frac{1}{j!}\phi^jB},
\quad
\frac{1}{n!}\phi^n({A}\wedge{B})=\sum_{k+j=n}({\frac{1}{k!}\phi^kA}\wedge {\frac{1}{j!}\phi^jB}).
\end{equation}
It is equivalent to
\begin{equation}\label{endo n}
\phi^n\lbra{A}{B}=\sum_{k=1}^nC^n_k\lbra{\phi^kA}{\phi^{n-k}B},
\quad
\phi^n({A}\wedge{B})=\sum_{k=1}^nC^n_k({\phi^kA}\wedge {\phi^{n-k}B}).
\end{equation}

Finally,  $\Phi(t)$ intertwines $\dbar_{\Gamma(t)}$ if and only if they satisfy the identity
(\ref{intertwine}). Assuming that $\Phi(t)$ is a homomorphism of the Courant bracket on $\wedge^\bullet L$,
 we need to show that for any section $A$ of $\wedge^\bullet L$,
\begin{equation}\label{intertwine with t}
\Phi(t)(\dbar A)+\lbra{\Phi(t)\Gamma(t)}{\Phi(t)A}=
\dbar(\Phi(t)A).
\end{equation}
Consider the infinite product $\Phi(t)\Gamma(t)$. Its degree-$n$ term is equal to
\[
\sum_{k+j=n}\frac{1}{j!}\phi^j(-1)^{k-1}\frac{1}{k!}\phi^{k-1}\Gamma_1
=\left(\sum_{k+j=n}(-1)^{k-1}\frac{1}{j!k!}\right)\phi^{n-1}\Gamma_1
\]
On the other hand, consider the power series
\[
g=\sum_{k=1}^\infty (-1)^k\frac{1}{k!}x^{k-1} \quad \mbox{ and } \quad
e^x=\sum_{k=1}^\infty\frac{1}{k!}x^{k}.
\]
We have
$1-xg=e^{-x}.$ Therefore, $e^x-xe^xg=1$. i.e. $e^x(xg)=e^x-1$. Equating the $n$-th
order terms for $n\geq 1$, we find that
\[
\sum_{k+j=n}(-1)^{k-1}\frac{1}{j!k!}=\frac{1}{n!}.
\]
Therefore, $\Phi(t)\Gamma(t)=\sum_{n\geq 1}\frac{1}{n!}t^n\phi^{n-1}\Gamma_1$.
Then the identity (\ref{intertwine with t}) at degree $n$ becomes
\[
\frac{1}{n!}\phi^n\dbar A+\sum_{k+j=n}\lbra{\frac{1}{k!}\phi^{k-1}\Gamma_1}{\frac{1}{j!} \phi^j A}
=\frac{1}{n!}\dbar \phi^n A
\]
Equivalently, it is
\begin{equation}\label{intertwine n}
\dbar \phi^n A-\phi^n\dbar A=\sum_{k=1}^nC^n_k\lbra{\phi^{k-1}\Gamma_1}{ \phi^{n-k} A}
\end{equation}

To complete a proof of Theorem \ref{thm:integrable}, we need to prove that the identities
(\ref{MC n}), (\ref{endo n}) and (\ref{intertwine n}) hold.

\begin{lemma} Suppose that $\Gamma_1$ and $\phi$ form a compatible pair,
then  identity {\rm (\ref{MC n})} holds.
\end{lemma}
\bproof
Using the fact that $\dbar\Gamma_1=0$ and equation
(\ref{infinitesimal 1}), we get a telescopic sum
\begin{eqnarray*}
\dbar\phi^{n-1}\Gamma_1
&=&\dbar(\phi^{n-1}\Gamma_1)-\phi\dbar\phi^{n-2}\Gamma_1\\
&+&\phi\dbar\phi^{n-2}\Gamma_1-\phi^2\dbar\phi^{n-3}\Gamma_1
+\phi^2\dbar\phi^{n-3}\Gamma_1-\phi^3\dbar\phi^{n-4}\Gamma_1\\
                       &+&\dots\ \ \dots \ \ \dots
                         +\phi^{n-2}\dbar\phi\Gamma_1-\phi^{n-1}\dbar\Gamma_1\\
&=&\sum_{h=0}^{n-2}\phi^h\lbra{\Gamma_1}{\phi^{n-2-h}\Gamma_1}.
\end{eqnarray*}
Since $\phi$ satisfies (\ref{phi endo}) and the Schouten bracket is
commutative when restricted to section of $\Lambda^2L$, we rewrite the above identity as
\begin{eqnarray*}
\dbar\phi^{n-1}\Gamma_1
&=&\sum_{h=0}^{n-2}\sum_{k=0}^{h}C^h_k\lbra{\phi^k\Gamma_1}{\phi^{n-2-k}\Gamma_1}
=\sum_{h=1}^{n-1}\sum_{k=1}^{h}C^{h-1}_{k-1}\lbra{\phi^{k-1}\Gamma_1}{\phi^{n-1-k}\Gamma_1}\\
  &=&\sum_{k=1}^{n-1}\left(\sum_{h=k}^{n-1}C^{h-1}_{k-1}\right)
  \lbra{\phi^{k-1}\Gamma_1}{\phi^{n-1-k}\Gamma_1}\\
  &=&\sum_{k=1}^{n-1}
  C^{n-1}_{k}\lbra{\phi^{k-1}\Gamma_1}{\phi^{n-1-k}\Gamma_1}.
\end{eqnarray*}
Performing the index substitution $k\mapsto n-k$ and using the
commutativity of the Schouten bracket again, we get
\begin{eqnarray*}\dbar\phi^{n-1}\Gamma_1&=&\frac{1}{2}\sum_{k=1}^{n-1}
\left(C^{n-1}_{k}+C^{n-1}_{k-1}\right)\lbra{\phi^{k-1}\Gamma_1}{\phi^{n-1-k}\Gamma_1}\\
&=&\frac{1}{2}\sum_{k=1}^{n-1}
C^n_k\lbra{\phi^{k-1}\Gamma_1}{\phi^{n-1-k}\Gamma_1}.
\end{eqnarray*}
\eproof

\begin{lemma}\label{endo to iso}
 Suppose that $\Gamma_1$ and $\phi$ form a compatible pair, then the two identities
 in {\rm (\ref{endo n})} hold.
\end{lemma}
\bproof
It is an elementary induction.
The proof for both cases are identical. We work only through the
case with Schouten bracket. When $n=1$, the equation (\ref{endo n}) is
precisely the equation (\ref{phi endo}), which is satisfied by assumption. Assuming that the equation (\ref{endo n}) holds for
all $k\leq n$. We next compute $\phi^{n+1}\lbra{A}{B}$, which we
take as $\phi(\phi^n\lbra{A}{B})$. By induction hypothesis, it is
equal to
\begin{eqnarray*}
&&\sum_{k=0}^nC_k^n\phi(\lbra{\phi^{n-k}A}{\phi^kB})\\
&=&\sum_{k=0}^nC_k^n(\lbra{\phi^{n+1-k}A}{\phi^kB}+\lbra{\phi^{n-k}A}{\phi^{k+1}B})\\
&=&\lbra{\phi^{n+1}A}{B}+\sum_{k=1}^nC_k^n\lbra{\phi^{n+1-k}A}{\phi^kB}\\
&&+\sum_{k=0}^{n-1}C_k^n\lbra{\phi^{n-k}A}{\phi^{k+1}B}+
\lbra{A}{\phi^{n+1}B}\\
&=&\lbra{\phi^{n+1}A}{B}+\sum_{k=1}^n(C_k^n+C_{k-1}^n)
\lbra{\phi^{n+1-k}A}{\phi^kB}+\lbra{A}{\phi^{n+1}B}\\
&=&\sum_{k=0}^{n+1}C_k^{n+1}\lbra{\phi^{n-k} A}{\phi^{k} B}.
\end{eqnarray*}
\eproof

\begin{lemma} Suppose that $\Gamma_1$ and $\phi$ form a compatible pair,
then  identity {\rm (\ref{intertwine n})} holds
for all $n\geq 1$.
\end{lemma}
\bproof
Since $\dbar\Gamma_1=0$, we substitute $A$ by $\Gamma_1$ in (\ref{infinitesimal 1}) to see that
 identity
(\ref{intertwine n}) holds when $n=1$. Assume that (\ref{intertwine n}) holds for
all $k\leq n$. We next prove that it holds for $n+1$. Since
\begin{eqnarray*}
&&\dbar\phi^{n+1}A-\phi^{n+1}\dbar A\\
&=&\dbar\phi^n(\phi A)-\phi^n(\dbar\phi A)+\phi^n\left((\dbar\phi
A)-\phi\dbar A\right),
\end{eqnarray*}
by  induction hypothesis, the above is equal to
\[
\sum_{k=1}^nC_k^n\lbra{\phi^{k-1}\Gamma_1}{\phi^{n+1-k}A}
+\phi^n\lbra{\Gamma_1}{A}.
\]
By Lemma \ref{endo to iso}, it is equal to
\begin{eqnarray*}
&&\sum_{k=1}^nC_k^n\lbra{\phi^{k-1}\Gamma_1}{\phi^{n+1-k}A}
+\sum_{k=0}^nC_k^n\lbra{\phi^k\Gamma_1}{\phi^{n-k}A}\\
&=&\sum_{k=1}^n\left(C_k^n+C_{k-1}^n\right)\lbra{\phi^{k-1}\Gamma_1}{\phi^{n+1-k}A}
+\lbra{\phi^n\Gamma_1}{A}.
\end{eqnarray*}
By Pascal Identity, it is equal to
\[
\sum_{k=1}^nC_k^{n+1}\lbra{\phi^{k-1}\Gamma_1}{\phi^{n+1-k}A}+\lbra{\phi^n\Gamma_1}{A}
=\sum_{k=1}^{n+1}C_k^{n+1}\lbra{\phi^{k-1}\Gamma_1}{\phi^{n+1-k}A}.
\]
\eproof

\section{Holomorphic Poisson manifolds}\label{sec:Poisson}

On a complex manifold $(M,J)$,
$
L=T^{1,0}\oplus T^{*(0,1)}$, and $\oL=L^*=T^{0,1}\oplus
T^{*(1,0)}.$
Therefore, the exterior bundle has a decomposition
\[
\wedge^\bullet L=\oplus_k\left(\oplus_{p+q=k}\wedge^pT^{1,0}\otimes \wedge^qT^{*(0,1)}  \right)
\]
We will  use the notations $T^{p,0}=\wedge^pT^{1,0}$ and
$T^{*(0,q)}=\wedge^qT^{*(0,1)}$. Sections of $T^{p,0}$ are addressed as $(p,0)$-vectors, more generally
polyvector fields.

\subsection{Type decomposition of deformations}
The  cohomology of $\DGA (J)$
decomposes accordingly into the direct sum of classical Dolbeault cohomology
with the  sheaf of exterior product of the holomorphic tangent
bundle as coefficients.
\begin{equation}
H^k_J=\oplus_{p+q=k, p,q\geq 0}H^q(M, T^{p,0}).
\end{equation}
If $\Gamma_1$ is in $H^2_J$ it has three components:
\begin{equation}
\Gamma_1=\Lambda+{\widehat\Gamma}_1+\Omega\in
H^0(M, T^{2,0})\oplus H^1(M, T^{1,0})\oplus
H^2(M, {\mathcal O}),
\end{equation}
where $\Lambda$ is a (2,0)-bivector field,
$\Omega$ is a (0,2)-form, and
${\widehat\Gamma}_1$ is a classical infinitesimal complex deformation.
Similarly,
\begin{eqnarray*}
L^*\otimes L &=&\END( L, L)\\
&=&\END (T^{1,0},T^{1,0}) \oplus \END(T^{*(0,1)},T^{*(0,1)})\\
&& \quad \quad \oplus \END(T^{1,0}, T^{*(0,1)})\oplus
\END(T^{*(0,1)}, T^{1,0}).
\end{eqnarray*}

If $\phi$ is a section of  $L^*\otimes L$, we represent its decomposition
by $\phi=\phi_1+\phi_2+\phi_3+\phi_4$ such that
\begin{eqnarray*}
&\phi_1\in C^\infty(\END (T^{1,0},T^{1,0})), \quad
\phi_2\in C^\infty(\END(T^{*(0,1)},T^{*(0,1)})),&\\
&\phi_3\in C^\infty(\END(T^{1,0}, T^{*(0,1)})), \quad \phi_4\in
C^\infty(\END(T^{*(0,1)}, T^{1,0})).&
\end{eqnarray*}

\begin{proposition} A pair $\Gamma_1\in C^\infty(M, \wedge^2L)$ and $\phi\in
C^\infty(M, L^*\otimes L)$ is compatible if and only if the pairs $(\Lambda,
\phi_4)$, $(\Omega, \phi_3)$ and $({\widehat\Gamma}_1, \phi_1+\phi_2)$ are
compatible.
\end{proposition}
\bproof This theorem is an inspection of  type
decompositions. For example,
\[
\dbar\Gamma_1=\dbar\Lambda+\dbar{\widehat\Gamma}_1+\dbar\Omega.
\]
Since
$\dbar\Lambda\in$ $C^\infty(M, T^{2,0}\otimes T^{*(0,1)})$,
$\dbar{\widehat\Gamma}_1$ $\in C^\infty(M, T^{1,0}\otimes T^{*(0,2)})$, and $\dbar\Omega\in$
 $ C^\infty(M,  T^{*(0,3)})$, each component has to vanish individually if
 $\dbar\Gamma_1=0.$ i.e.
\[
\dbar\Lambda=0, \quad \dbar{\widehat\Gamma}_1=0, \quad \dbar\Omega=0.
\]

Next, for all $Z\in C^\infty(M, T^{1,0})$ and $\oom\in C^\infty(M, T^{*(0,1)})$,
\begin{eqnarray*}
&\lbra{\Lambda}{Z}\in C^\infty(T^{2,0}), \quad \lbra{\Lambda}{\oom}\in
C^\infty(T^{(1,0)}\otimes T^{*(0,1)});&\\
&\lbra{{\widehat\Gamma}_1}{Z}\in C^\infty(T^{*(0,2)}), \quad
\lbra{{\widehat\Gamma}_1}{\oom}=0;&\\
 &\lbra{\Omega}{Z} \in C^\infty(T^{1,0}\otimes
T^{*(0,1)}), \quad \lbra{\Omega}{\oom} \in C^\infty(T^{*(0,2)}).&
\end{eqnarray*}
On the other hand,
\begin{eqnarray*}
&\dbar(\phi_1(Z))-\phi_1(\dbar Z)\in C^\infty(T^{1,0}\otimes T^{*(0,1)}),
\quad \dbar(\phi_1(\oom))-\phi_1(\dbar \oom)=0,&\\
&\dbar(\phi_2(Z))-\phi_2(\dbar Z)=-\phi_2(\dbar Z)\in
C^\infty(T^{1,0}\otimes T^{*(0,1)}), \quad \dbar(\phi_2(\oom))-\phi_2(\dbar
\oom)\in C^\infty(T^{*(0,2)}),&\\
&\dbar(\phi_3(Z))-\phi_3(\dbar Z)\in C^\infty(T^{*(0,2)}), \quad
\dbar(\phi_3(\oom))-\phi_3(\dbar \oom)=0,&\\
&\dbar(\phi_4(Z))-\phi_4(\dbar Z)=-\phi_4(\dbar Z)\in
C^\infty(T^{2,0}), \quad \dbar(\phi_4(\oom))-\phi_4(\dbar \oom)\in C^\infty(T^{1,0}\otimes
T^{*(0,1)}).&
\end{eqnarray*}
By equating the types, we arrive at the conclusion of this proposition.
\eproof

In view of the last proposition and the decomposition of $H^2_J$, one should
focus an initial analysis of deformations on
the simple types, namely those whose infinitesimal deformations are
contained in a unique summand of the decomposition
of $H^2_J$.

Infinitesimal deformations given by a $\dbar$-closed section ${\widehat\Gamma}_1$
of $T^{1,0}\otimes T^{*(0,1)}$ could always be
represented and analyzed as classical complex deformation theory.

If one considers a $\dbar$-closed 2-form representing an element in  $H^2(M, \calo)$,
then by definition of Courant bracket $\lbra{\Omega}{\Omega}=0$.
Therefore, $\Omega$ satisfies the Maurer-Cartan equation, and
\[
L_{\overline\Omega}=\{X+\alpha+\iota_X{\overline\Omega}: X+\alpha\in L\}
\]
is a generalized complex structure. The issue of integrability is trivial.
 However, this deformation does not change the type of the generalized complex
structure.
It is still type-$n$ where $n$ is the complex dimension of the manifold $M$.
If the $(0,2)$-form $\Omega$ is not only $\dbar$-closed but also closed, then
this deformation is trivial within the realm of
generalized complex structures because the deformation is only the result of a B-field transformation
\cite{Marco}.

\subsection{Holomorphic bivector fields}
Suppose that $\Gamma$ is a deformation whose first order term is a
bivector field $\Lambda$ with $\dbar\Lambda=0$.
Let $\Gamma_2$ be its second order term. As $\Gamma$ satisfies
the Maurer-Cartan equation, up to second order term, we have
\[
\dbar(t\Lambda+t^2\Gamma_2)+\frac12\lbra{t\Lambda+t^2\Gamma_2}{t\Lambda+t^2\Gamma_2}=0.
\]
It yields
\[
\dbar\Gamma_2+\frac12\lbra{\Lambda}{\Lambda}=0.
\]
Since $\Lambda$ is a bivector, $\lbra{\Lambda}{\Lambda}$ is a $(3,0)$-vector field.
On the other hand, $\dbar\Gamma_2$ must have a components with
$(0,1)$-forms. Therefore, the only solution is when $\lbra{\Lambda}{\Lambda}=0$.
It follows immediately that $\Lambda$ is a solution of the
Maurer-Cartan equation and $\Gamma_2$ could be chosen to be zero.
Therefore, a bivector field $\Lambda$ representing an element in $H^0(M, T^{2,0})$ is an
infinitesimal deformation of an integrable deformation if and only if $\lbra{\Lambda}{\Lambda}=0$.

\begin{definition} A $(2,0)$-vector field is a holomorphic Poisson structure on a
complex manifold if $\dbar\Lambda=0$ and
$\lbra{\Lambda}{\Lambda}=0$. In such case, we call $\Lambda$ a holomorphic Poisson
vector field.
\end{definition}

Given such a $\Lambda$, suppose that $\phi\in C^\infty(\END(T^{*(0,1)}, T^{1,0}))$ is
compatible with $\Lambda$. By Theorem \ref{thm:integrable},
$\Phi=\sum \frac{1}{n!}\phi^n$ is a $\DGA$ homomorphism. However, as an endomorphism
from the bundle $L$ to $L$, its kernel contains at least
$T^{1,0}$. Therefore, $\phi\Lambda=0$ and $\phi\circ\phi=0$. Therefore, we could
conclude that the homomorphism $\Phi$ is simply $1+\phi$.
Furthermore, given a section $X+\overline\alpha$ of
$T^{1,0}\oplus T^{*(0,1)}$, $\Phi(X+{\overline\alpha})=X+\phi({\overline\alpha})+{\overline\alpha}$.
As the vector part is $X+\phi({\overline\alpha})$ and the form part is
${\overline\alpha}$, $X+{\overline\alpha}$ is in the kernel of
$\Phi$ if and only if it is identically zero. Therefore, $\Phi$ as a
bundle map from $L$ to $L$ is an isomorphism. It is extended to an
isomorphism from the exterior bundle $\wedge^\bullet L$ to $\wedge^\bullet L$.
Therefore, $\Phi$ is not only a $\DGA$ homomorphism, but also an isomorphism.
We summarize our observation below.

\begin{theorem}\label{thm:poisson}
Let $M$ be a complex manifold with a holomorphic Poisson vector field $\Lambda$.
Suppose that $\phi$ is a section of $\END(T^{*(0,1)}, T^{1,0})$ compatible
with the $\Lambda$ in the sense of Definition {\rm \ref{compatible pair}}. Then $\Lambda$ defines a family of
generalized complex deformation
 of the complex structure on $M$ with $t\Lambda$. Moreover, if $\DGA(t\Lambda)$ represents
 the $\DGA$ of the deformed complex structure, then they are all isomorphic to  $\DGA(0)$,
 the differential Gerstenhaber algebra of the complex structure on the manifold $M$.
 \end{theorem}

 Although from the viewpoint of deformation of $\DGA$s, the presence of a compatible pair
 on a holomorphic Poisson manifold makes the deformation of $\DGA$s trivial, on the geometric
 level,  it is non-trivial. Recall that
 \[
 L_{\overline\Lambda}=\{X+{\overline\alpha}+\iota_{\overline\alpha}{\overline\Lambda}:
  X+{\overline\alpha}\in T^{1,0}\oplus T^{*(0,1)}\}.
 \]
 As $\iota_{\overline\alpha}{\overline\Lambda}$ is a (0,1)-vector, the type of
 the generalized complex structure
 $L_{\overline\Lambda}$ is different from the un-deformed one $L$.
 If $\Lambda$ as a bundle map from $T^{*(0,1)}$ to
 $T^{1,0}$ is everywhere non-degenerate, then $L_{\overline\Lambda}$ is a
 type-0 generalized complex structure. By a Gualtieri's lemma
 \cite{Marco}, there exists a symplectic structure $\Omega$ on the manifold $M$
 such that the complexified $\DGA$ of $\omega$ is isomorphic to that
 of $\DGA(\Lambda)$ via a B-field transformation. Since $\DGA(\Lambda)$ is isomorphic
 to $\DGA(0)$. We obtain the following result.

 \begin{theorem}\label{thm:weak mirror}
 Let $M$ be a manifold with complex  structure $J$. Denote its associated $\DGA$
  by $\DGA(J)$. Suppose that $\Lambda$ is a non-degenerate holomorphic Poisson
 structure. If there exists  a section of $T^{1,0}\otimes T^{0,1}$ compatible with
  $\Lambda$ in the sense of Definition {\rm \ref{compatible pair}}, then there exists a
 symplectic structure $\Omega$ in the deformation family of $J$ such that $\DGA(\Omega)$ is isomorphic to
 $\DGA(J)$.
 \end{theorem}

 In the sense of Merkulov, the pair $(M, J)$ and $(M, \Omega)$ form a weak mirror pair
  \cite{CLP} \cite{COP} \cite{Mer-note}.

\subsection{Rational surfaces}
In this section, we compute the first cohomology of some well known holomorphic Poisson manifolds to demonstrate
that for many holomorphic Poisson structures, Theorem \ref{thm:poisson} does not have solution.

Assume that we have a compact holomorphic Poisson manifold. Denote the Poisson bivector field by $\Lambda$.
Consider  $Z$ a section of $T^{1,0}$ and $\oomega$ a section of $T^{*(0,1)}$. Then $Z+\oomega$ is a section of
$L=T^{1,0}\oplus T^{*(0,1)}$. By Proposition \ref{reference}, it represents an element of the first cohomology
of $\DGA(\dbar_{\Lambda})$ if and only if
$\dbar_{\Lambda}(Z+\oomega)=0$. That is
\begin{eqnarray*}
\dbar_{\Lambda}(Z+\oomega)&=&\dbar Z+\lbra{\Lambda}{Z}+\dbar\oomega+\lbra{\Lambda}{\oomega}\\
&=&\lbra{\Lambda}{Z}+\dbar Z+\lbra{\Lambda}{\oomega}+\dbar\oomega=0.
\end{eqnarray*}
The terms above are sections of $\wedge^2L=T^{2,0}\oplus T^{1,0}\otimes T^{*(0,1)}\oplus T^{*(0,2)}$.
As each component in this decomposition has to vanish, we conclude that
\begin{equation}\label{dbar closed}
\lbra{\Lambda}{Z}=0, \quad \dbar Z+\lbra{\Lambda}{\oomega}=0, \quad \dbar\oomega=0.
\end{equation}
In particular the (0,1)-form $\oomega$ is $\dbar$-closed. To push this computation further,
assume that the Dolbeault cohomology $H^1(M, \calo)$ vanishes.
It follows that the (0,1)-form is $\dbar$-exact, and there is a smooth function $f$ on the manifold
$M$ such that $\oomega=\dbar f$. Consider the vector field $V=\lbra{\Lambda}{f}$. Since $\dbar\Lambda=0$,
\[
\dbar V =\dbar \lbra{\Lambda}{f}=\lbra{\dbar\Lambda}{f}-\lbra{\Lambda}{\dbar f}
=-\lbra{\Lambda}{\oomega}.
\]
With (\ref{dbar closed}) above, we conclude that $\dbar (Z-V)=0$. Therefore, $Z-V$ is a holomorphic vector
field on the manifold $M$.
 By Jacobi identity of Gerstenhaber algebras,
\[
\lbra{\Lambda}{\lbra{\Lambda}{f}}+\lbra{\Lambda}{\lbra{f}{\Lambda}}+\lbra{f}{\lbra{\Lambda}{\Lambda}}=0.
\]
Since $\lbra{\Lambda}{\Lambda}=0$, the above is reduced to
\[
\lbra{\Lambda}{V}=\lbra{\Lambda}{\lbra{\Lambda}{f}}=0.
\]
Combined with the first identity in (\ref{dbar closed}), we conclude that
\[
\lbra{\Lambda}{Z-V}=0.
\]
Let $W=Z-V$, then $Z=W+V=W+\lbra{\Lambda}{f}$ such that $\dbar W=0$ and $\lbra{\Lambda}{W}=0$.
Moreover, the section
\[
Z+\oomega=W+\lbra{\Lambda}{f}+\dbar f=W+\dbar_{\Lambda}f.
\]
Since $\dbar_{\Lambda}f$ is $\dbar_{\Lambda}$-exact, $W$ and $Z+\oomega$ represent the same cohomology
class in $H^1_{\dbar_{\Lambda}}$.

\begin{proposition}\label{first cohomology} Suppose that $M$ is a holomorphic Poisson manifold with
Poisson vector field $\Lambda$. If $H^1(M, \calo)$ vanishes, then
\[
H^1_{\dbar_{\Lambda}}=\{W\in H^0(M, T^{1,0}): \lbra{\Lambda}{W}=0\}.
\]
\end{proposition}

On the other hand, the first cohomology of $\DGA(\dbar)$ is equal to
\[
H^0(M, T^{1,0})\oplus H^1(M, \calo).
\]
Given the assumption of Proposition \ref{first cohomology}, it is equal to $H^0(M, T^{1,0})$.
It is easy to find example on which there exists non-trivial holomorphic Poisson structures but it does
not admit compatible pairs due to the
difference between $H^0(M, T^{1,0})$ and $H^1_{\dbar_{\Lambda}}$.
For instance, there is a classification of compact complex surfaces admitting holomorphic Poisson structures \cite{BM}.
Among them, the minimal rational surfaces are all holomorphic Poisson manifolds with vanishing irregularity.
Except when the surface is a complex projective plane, they are rational ruled surfaces.

For the complex projective plane $\Lambda$ is an element in $H^0(\cp^2, \calo(3))$. It could be identified
to a homogeneous polynomial of degree-3 in the homogeneous coordinates  of the
complex projective plane. Meanwhile the space of holomorphic vector fields $H^0(\cp^2, T^{1,0})$
is the complex algebra $\mathfrak{sl}(3, \CC)$, treated as the set of $3\times 3$-matrices acting on
of $\CC^3$ by natural matrix multiplications. From this perspective, for any $W$ in
$H^0(\cp^2, T^{1,0})$,  the action $\lbra{W}{-}$ on $H^0(\cp^2, \calo(3))$
is the induced representation of $\mathfrak{sl}(3, \CC)$ on the third symmetric product $S^3\CC^3$.
Then for each $\Lambda\neq 0$, one could find a $W$ such that
  $\lbra{W}{\Lambda}\neq 0$.
Therefore, for each holomorphic Poisson structure on the complex projective plane, $H^1_{\dbar_{\Lambda}}$ is strictly smaller then
 $H^1_{\dbar}=H^0(\cp^2, T^{1,0})=\mathfrak{sl}(3, \CC)$. It shows that  $\DGA(\cp^2, \dbar)$ and
$\DGA(\cp^2, \dbar_{\Lambda})$ for any holomorphic Poisson structure could never be quasi-isomorphic.

\subsection{Hopf surfaces}
In this section, we compute  $H^1_{\dbar_{\Lambda}}$ when the underlying manifold
$M$ is the Hopf surface, and demonstrates
that this does admit compatible pairs.

Consider $\CC^2$ with coordinates $z=(z_1, z_2)$. Let $\lambda>1$ be a real number. It generates a one-parameter
group of automorphism on $\CC^2$. The quotient of $\CC^2\backslash \{ 0\}$ with respect this group is diffeomorphic
to the Lie group $M=U(1)\times SU(2)$. The complex structure on $\CC^2$ descends onto $M$ to define an
integrable complex structure, invariant of the left-action of the Lie group. In this section, by Hopf surface, we
mean this particular complex structure. The classical complex deformation theory of this
complex structure was analyzed by Dabrowski
\cite{Dab}. We focus on the deformations generated by its holomorphic Poisson structures.
Consider
\begin{eqnarray*}
X_0 = \frac12(z_1\frac{\partial}{\partial z_1}+z_2\frac{\partial}{\partial z_2}), &\quad&
X_1 = \frac{i}2(z_1\frac{\partial}{\partial z_1}-z_2\frac{\partial}{\partial z_2})\\
X_2 = \frac{i}2(z_2\frac{\partial}{\partial z_1}+z_1\frac{\partial}{\partial z_2}), &\quad&
X_3 = \frac12(-z_2\frac{\partial}{\partial z_1}+z_1\frac{\partial}{\partial z_2}),
\end{eqnarray*}
and
\[
\overline{\sigma} = \overline{\partial} \ln |z|^2 = \frac{z_1d{\overline z}_1+z_2d{\overline z}_2}
{|z_1|^2+|z_2|^2}.
\]
The cohomology spaces for the $\DGA(J)$ are given below. The computation of these cohomology spaces
are not new. We do not present any details.
\begin{eqnarray}
&H^1(M, \co)=\langle {\overline\sigma} \rangle, \quad
H^0(M, T^{1,0})=\langle X_0, X_1, X_2, X_3\rangle\cong {\mathfrak{u}}(1)\oplus
{\mathfrak{sl}}(2), & \label{hopf h1} \\
&H^1(M, T^{1,0})=\langle X_0\wedge{\overline\sigma}, X_1\wedge{\overline\sigma},
X_2\wedge{\overline\sigma}, X_3\wedge{\overline\sigma}\rangle, &\\
&H^0(M, T^{2,0})=\langle X_0\wedge X_1, X_0\wedge X_2, X_0\wedge X_3\rangle,&
\label{hopf poisson}
\\
&H^1(M, T^{2,0})=\langle X_0\wedge X_1\wedge{\overline\sigma},
X_0\wedge X_2\wedge{\overline\sigma}, X_0\wedge X_3\wedge{\overline\sigma} \rangle.&
\end{eqnarray}
In addition,
\begin{eqnarray}
&\lbra{X_0}{X_1}=0, \quad \lbra{X_0}{X_2}=0, \quad \lbra{X_0}{X_3}=0,&\label{u1 su2}\\
&\lbra{X_1}{X_2}=-X_3, \quad \lbra{X_2}{X_3}=-X_1, \quad \lbra{X_3}{X_1}=-X_2.&
\end{eqnarray}
Set $f= \ln |z|^2$, then $\cl_{X_0}f=\frac12$. For $j=1,2,3$, define $f_j=\cl_{X_j}f$, then
\[
f_1= \frac{i}{2|z|^2}(z_1\zbar_1-z_2\zbar_2), \quad
f_2=\frac{i}{2|z|^2}(z_2\zbar_1+z_1\zbar_2), \quad
f_3=\frac{1}{2|z|^2}(-z_2\zbar_1+z_1\zbar_2).
\]
The functions $f_1, f_2, f_3$ are invariant of the group of actions generated by
$(\lambda z_1, \lambda z_2)$, and hence they are globally defined on the quotient space $M$.
Then we have
\begin{equation}\label{X osigma}
\lbra{X_0}{\osigma}=0, \quad
\lbra{X_1}{\osigma}=\dbar f_1, \quad
\lbra{X_2}{\osigma}=\dbar f_2, \quad
\lbra{X_3}{\osigma}=\dbar f_3.
\end{equation}
Whenever $A=a_1X_1+a_2X_2+a_3X_3$ is a holomorphic vector field in the
${\mathfrak{sl}}(2)$ component of $H^0(M, T^{1,0})$,
\[
\lbra{A}{\osigma}= a_1\lbra{X_1}{\osigma}+a_2\lbra{X_2}{\osigma}+a_3\lbra{X_3}{\osigma}
=\dbar(a_1f_1+a_2f_2+a_3f_3).
\]
We use the notation $f_A$ to denote the function $a_1f_1+a_2f_2+a_3f_3$. By (\ref{X osigma}),
\begin{equation}\label{A osigma}
\lbra{A}{\osigma}=\dbar f_A.
\end{equation}
Since $X_0$ commutes with $X_j$ for $j=1,2,3$,
$
\cl_{X_0}f_j=\cl_{X_j}\cl_{X_0}f=\cl_{X_j}\frac12=0.
$
Then for all $A$
\begin{equation}\label{X_0 f_A}
\cl_{X_0}f_A=0.
\end{equation}

Given the above preparation, we begin to compute the first cohomology of $\DGA(\dbar_{\Lambda})$ where $\Lambda$ is any holomorphic Poisson structure on $M$.
Let $A=a_1X_1+a_2X_2+a_3X_3$ be a holomorphic vector field, then $\Lambda=X_0\wedge A$ is a holomorphic
Poisson structure. As noted in (\ref{hopf poisson}),
by choosing the complex numbers $(a_1, a_2, a_3)$, we exhaust all holomorphic Poisson structure.

Now we calculate the first cohomology with respect to $\dbar_{\Lambda}=\dbar+\lbra{\Lambda}{-}$.
Suppose that $Z$ is a smooth (1,0)-vector field and $\oomega$ is a smooth (0,1)-form.
$\dbar_{\Lambda}(Z+\oomega)=0$ if and only if $Z$ and $\oomega$ satisfy the constraints
(\ref{dbar closed}). Once again, they are
\begin{equation}\label{three conditions}
\dbar\oomega=0, \quad \lbra{\Lambda}{Z}=0, \quad
\dbar Z+\lbra{\Lambda}{\oomega}=0.
\end{equation}
Since the cohomology $H^1(X, \calo)$ is spanned by $\osigma$, there exists a function $\psi$ and a
unique complex
number $a$ such that
\[
\oomega=a\osigma+\dbar\psi.
\]
Let $V$ be the vector field $\lbra{\Lambda}{\psi}$. Since $\dbar\Lambda=0$,
\begin{eqnarray*}
\dbar V&=&\dbar \lbra{\Lambda}{\psi}= \lbra{\dbar\Lambda}{\psi}-\lbra{\Lambda}{\dbar\psi}
=-\lbra{\Lambda}{\dbar\psi}\\
&=&-\lbra{\Lambda}{\oomega-a\osigma}=-\lbra{\Lambda}{\oomega}+a\lbra{\Lambda}{\osigma}.
\end{eqnarray*}
By definition of $\Lambda$, (\ref{X osigma}) and (\ref{A osigma}), this is equal to
\[
-\lbra{\Lambda}{\oomega}+aX_0\wedge\lbra{A}{\osigma}
=-\lbra{\Lambda}{\oomega}+aX_0\wedge \dbar f_A
=-\lbra{\Lambda}{\oomega}-a\dbar(f_A X_0).
\]
It follows from (\ref{three conditions}) that
\begin{equation}
\dbar (V-Z+a f_AX_0)=0.
\end{equation}
Next, consider the Schouten bracket. By (\ref{three conditions}),
\begin{eqnarray*}
&&\lbra{\Lambda}{V-Z+a f_AX_0}\\
&=&\lbra{\Lambda}{\lbra{\Lambda}{\psi}}-\lbra{\Lambda}{Z}+\lbra{X_0\wedge A}{af_A X_0}\\
&=&\lbra{\Lambda}{\lbra{\Lambda}{\psi}}+aX_0\wedge \lbra{A}{f_A X_0}-aA\wedge \lbra{X_0}{f_A X_0}\\
&=&\lbra{\Lambda}{\lbra{\Lambda}{\psi}}+aX_0\wedge \lbra{A}{f_A} X_0
+a f_AX_0\wedge \lbra{A}{X_0}-aA\wedge \lbra{X_0}{f_A} X_0.
\end{eqnarray*}
Due to (\ref{u1 su2}) and (\ref{X_0 f_A}), this is equal to
$\lbra{\Lambda}{\lbra{\Lambda}{\psi}}.$
By the Jacobi identity for Gerstenhaber algebra and the fact that $\lbra{\Lambda}{\Lambda}=0$,
 $\lbra{\Lambda}{V-Z+a f_AX_0}=0$. Define
\begin{equation}
W=-V+Z-a f_AX_0=-\lbra{X_0\wedge A}{\psi}+Z-af_AX_0.
\end{equation}

Then $\lbra{\Lambda}{W}=0$. However, by (\ref{three conditions}) and the identity above,
\begin{equation}
\lbra{\Lambda}{W}=-\lbra{X_0\wedge A}{W}=-X_0\wedge\lbra{A}{W}.
\end{equation}
As $\dbar W=0$, it is a linear combination of $X_0, X_1, X_2, X_3$.
Therefore, $\lbra{\Lambda}{W}$
is equal to zero if and only if there exist constants $b$ and $c$ such that
$W=bX_0+cA.$
Therefore,
\begin{equation}
Z=V+W+af_AX_0=\lbra{\Lambda}{\psi}+bX_0+cA+af_AX_0.
\end{equation}
As we have already resolved the first
two constraints in (\ref{three conditions}), we could now substitute
$Z$ in the last constraint to check that it does not generate
additional conditions.
So, $Z+\oomega$ is $\dbar_{\Lambda}$-closed for $\Lambda=X_0\wedge A$
if and only if there exist a function $\psi$ and
constants $a, b, c$ such that
\begin{eqnarray*}
\oomega &=&a\osigma+\dbar\psi,\\
Z &=&\lbra{\Lambda}{\psi}+bX_0+cA+af_AX_0=\lbra{\Lambda}{\psi}+W+af_AX_0.
\end{eqnarray*}
Since
\[
Z+\oomega=\lbra{\Lambda}{\psi}+bX_0+cA+af_AX_0+a\osigma+\dbar\psi
=bX_0+cA+af_AX_0+a\osigma+\dbar_{\Lambda}\psi,
\]
 $Z+\oomega$ and
$bX_0+cA+af_AX_0+a\osigma$ represent the same element in the first cohomology
space $H^1(M, \dbar_{X_0\wedge A})$.
Therefore, we have
\begin{equation}
H^1(X, \dbar_{X_0\wedge A})=\langle X_0, A, f_AX_0+\osigma\rangle\cong \CC^3.
\end{equation}
On the other hand, it is noted in (\ref{hopf h1}) that the first cohomology of
$\DGA(J)$ is a five-dimensional space.
\[
H^1(M, \dbar)=
H^1(M, \co)\oplus H^0(M, T^{1,0})=\langle {\overline\sigma} \rangle
\oplus \langle X_0, X_1, X_2, X_3\rangle.
\]
Therefore, along the deformation given by holomorphic Poisson vector field
$\Lambda=X_0\wedge A$, the first cohomology jumps and hence $\Lambda$ could not be
part of any compatible pair.

\section{Holomorphic symplectic algebras}\label{sec:algebra}

In an explicit computation in \cite{Poon}, part of the result in Theorem \ref{thm:weak mirror} has
been observed on the Kodaira-Thurston surface. It was possible to do an explicit computation
due to the fact that  the manifold is a low-dimension nilmanifold.

If $H$ is a simply-connected nilpotent Lie group and $K$ is a co-compact subgroup, then the
quotient manifold $M=H/K$ is said to be a nilmanifold.
Let $\lie h$ be the Lie algebra of the group $H$,
the Chevalley-Eilenberg differential $d$ determines a complex
\[
d:\wedge^k{\lie h}^*\to \wedge^{k+1}{\lie h}^*.
\]
It is known for a long time that the inclusion $\lie h^*$ as invariant section of
$T^*$ induces an isomorphism on the cohomology level
 \cite{Nomizu}. If the nilmanifold has an invariant symplectic structure $\Omega$,
one could therefore consider this inclusion as a quasi-isomorphism from the differential
 Gerstenhaber algebra with invariant objects $\DGA(\lie h, \Omega)$ to the manifold level
$\DGA(H/K, \Omega)$.

There were a series of attempt to attain a similar result for Dolbeault cohomology
\cite{CF}  \cite{CFP} \cite{CFGAU}  \cite{Rolle}. This body of research generates a
 collection of examples of nilmanifolds for which the inclusion of invariant sections
in the space of sections of the bundle $L=T^{1,0}\oplus T^{*(0,1)}$ induces a quasi-isomorphism
of $\DGA$s. Kodaira-Thurston surfaces is a prominent example with small
dimension. To illustrate the theory of the past few chapters, we now focus on $\DGA(\lie h, J)$
for some Lie algebra $\lie h$.

In our subsequent computation, we do not restrict $\lie h$ to being nilpotent, but will
construct algebras on which there is a good collection of geometric objects as in \cite{COP}.

\subsection{Pseudo-K\"ahler structures}
Let $(\mathfrak g, \omega)$ denote a real Lie algebra equipped
with a symplectic structure $\omega$. Let $V$ denote the underlying vector
space of $\mathfrak g$. We seek a linear map $\gamma: \mathfrak g\to
{\END}(V)$ such that for all $x,y, z\in \mathfrak g$,
\begin{eqnarray}
&& \gamma(x) y-\gamma(y) x= [x,y]; \label{t free}\\
&&  \omega(\gamma(x) y, z)+\omega(y, \gamma(x)
z)=0; \label{sym}\\
&& \gamma([x,y])=\gamma (x) \gamma(y)-\gamma(y)
\gamma(x). \label{flat}
\end{eqnarray}
The last condition requires $\gamma$ to be a representation. The
second condition means that it is a
symplectic.

If one uses $\gamma$ as a operator of vector fields on the Lie group of
the algebra $\lie g$, the last condition is equivalent to require $\gamma$ to be a flat connection.
Condition in (\ref{sym}) is to require the connection to be symplectic. The condition in (\ref{t free})
is to require the connection to be torsion-free.

Given the representation $\gamma$, one obtains a semi-direct product Lie
algebra $\lie h:=\lie g\ltimes V$ with a Lie bracket  defined by
\begin{equation}
\lbra{(x,0)}{(y,0)}=([x,y],0)  \quad
\lbra{(x,0)}{(0,v)}=(0, \gamma(x)v),
\end{equation}
for all $x,y\in \lie g$ and $v\in V$. Here we denote a generic element in $\lie g\ltimes V$ in terms
of the decomposition $(x,u)\in \lie g\oplus V$.

On the semi-direct product, consider the linear map.
\begin{equation}
J(x,y)=(-y,x).
\end{equation}
This is an almost complex structure. The $(1,0)$ vectors are given by
\begin{equation}
\lie h^{1,0}=\{(x,-ix)\in (\lie g\oplus V)_\CC: x\in \lie g\}.
\end{equation}
$J$ is an integrable complex structure due to (\ref{t free}) because
\[
\lbra{x-iJx}{y-iJy}=\lbra{(x, -ix)}{(y, -iy)}=([x,y],-i(\gamma(x)y-\gamma(y)x)).
\]
The symplectic structure $\omega$ induces three different symplectic
forms on the semi-direct product $\lie h$.
\begin{eqnarray}
&\Omega_1((x,u), (y,v)):=-\omega(x,v)-\omega(u,y),&\\
&\Omega_2((x,u), (y,v)):=\omega(x,y)-\omega(u,v),&\\
&\Omega_3((x,u), (y,v)):=\omega(x,y)+\omega(u,v).&
\end{eqnarray}
With respect to the complex structure $J$, $\Omega_c=\Omega_1+i\Omega_2$ is a closed
 (2,0)-form. It is non-degenerate in the sense that the contraction map
  \[
  V\mapsto \Omega_c(V, {\ }), \quad \Omega_c: \lie h^{1,0}\to \lie h^{*(1,0)}
   \]
   is non-degenerate.
   The pair $(\Omega_c, J)$ is called a complex symplectic structure on the algebra $\lie h$.
   Let $\Lambda$ be the inverse mapping of $\Omega_c$.
   \[
   \Lambda: \lie h^{*(1,0)}\to \lie h^{1,0}.
   \]
   It is a matter of definition that $\Lambda\in \wedge^2\lie h^{1,0}=\lie h^{2,0}$.
   Therefore, it could play the role of an invariant holomorphic Poisson structure.
   Indeed we have the following
   \begin{lemma} Let $\Lambda$ be the inverse of $\Omega_c$, then it satisfies the following.
   \begin{itemize}
\item For any $\alpha,\beta\in \lie h^{*(1,0)}$,
$
\Lambda(\alpha, \beta)=-\Omega_c(\Omega_c^{-1}(\alpha),
\Omega_c^{-1}(\beta)).
$
\item $\lbra{\Lambda}{\Lambda}=0.$
\item $\dbar\Lambda=0$.
\end{itemize}
\end{lemma}
\bproof Beyond tracing definitions, the first identity is an elementary application of the
algebraic properties of Gerstenhaber algebra. The second identity is equivalent to $d\Omega_c=0$.
The last is another application of the algebraic properties of Gerstenhaber algebra combined with
a type decomposition argument.
\eproof

The last lemma leads to the next.
\begin{lemma} Given a symplectic algebra $(\lie g, \omega)$ with a flat torsion-free symplectic
connection on the underlying vector space $V$ of $\lie g$, then the semi-direct product
$\lie h=\lie g\ltimes V$  has a holomorphic Poisson structure $(J, \Lambda=\Omega_c^{-1})$.
\end{lemma}

Given the above holomorphic Poisson structure, we consider the generalized deformation generated by
the holomorphic Poisson vector field $\Lambda$. It yields
\begin{equation}
L_{\overline\Lambda}=\lie h^{1,0}\oplus \{\overline\zeta+\overline\Lambda { } \overline\zeta:
\overline\zeta\in \lie h^{*(0,1)}\}.
\end{equation}
Since $\overline\Lambda: \lie  h^{*(0,1)}\to \lie h^{1,0}$ is an isomorphism with ${\overline\Omega}_c$
as its inverse,
\begin{eqnarray*}
L_{\overline\Lambda} &=&\lie h^{1,0}\oplus
\{ {\overline\Omega}_c(\overline{Y})+\overline{Y}:
\overline{Y}\in \lie h^{(0,1)}\}.
\end{eqnarray*}
Since ${\overline\Omega}_c$ is a (0,2)-form, for any (1,0)-vector $X$, ${\overline\Omega}_c(X)=0$.
Therefore, the above is equal to
\begin{eqnarray*}
&=&
\{ X+{\overline\Omega}_c(X)+\overline{Y}+{\overline\Omega}_c(\overline{Y}):
X\in \lie h^{1,0}, \overline{Y}\in \lie h^{0,1}\}\\
&=&
\{ V+\Omega_c(V):
V\in \lie h^{1,0}\oplus\lie h^{0,1}\}
=
\{ V+\Omega_1(V)-i\Omega_2(V): V\in \lie h_c\}\\
&=& e^{\Omega_1}\{ V-i\Omega_2(V): V\in \lie h_c\}.
\end{eqnarray*}
The last equality means that the deformed generalized complex structure
$L_{\overline\Lambda}$ is the $B$-field transformation by the closed 2-form $\Omega_1$
of the generalized complex structure defined by the symplectic form $\Omega_2$.
In conclusion, we have

\begin{proposition}\label{prop:symplectic 2}
 Given a symplectic algebra $(\lie g, \omega)$ with a flat torsion-free symplectic
connection on the underlying vector space $V$ of $\lie g$, then up to the $B$-field transformation with respect
to the closed 2-form $\Omega_1$, the generalized deformation of the
 classical complex structure by holomorphic Poisson structure $\Lambda=(\Omega_1+i\Omega_2)^{-1}$ is the
 the symplectic structure $\Omega_2$. In particular,
 $\DGA(L_{\overline\Lambda})$ is isomorphic to $\DGA(\Omega_2)$.
 \end{proposition}

\subsection{Compatible pairs}
A different perspective in understanding $\DGA(L_{\overline\Lambda})$ is in terms of compatible pair.
That is to identify an element $\phi$ in $\lie h^{0,1}\otimes \lie h^{1,0}$
so that $(\Lambda, \phi)$ forms a compatible pair.

As $\Omega_3$ is a (1,1)-form and its contraction map is non-degenerate
   \[
   \Omega_3:\lie h^{1,0}\to \lie h^{*(0,1)},
   \]
   its inverse map
   \[
   \Omega_3^{-1}: \lie h^{*(0,1)}\to \lie h^{1,0}
   \]
   is a natural candidate to form a compatible pair with $\Lambda$.

   On the other hand, if $g$ is a non-degenerate symmetric bilinear
    form the algebra $\lie g$, it induces a non-degenerate form on $\lie g\ltimes V$ by
   \[
   \Delta((x,u), (y,v))=g(x,y)+g(u,v).
   \]
   Then its fundamental form is a (1,1)-form:
   \[
   \Omega_4((x,u), (y,v))=\Delta(J(x,u), (y,v))=\Delta((-u,x), (y,v))=g(x,v)-g(y,u).
   \]
   Therefore, $\Omega_4^{-1}$ is also a candidate to match with $\Lambda$ as a compatible pair.
  It is a natural question to ask when $\Omega_4$ is closed. It amounts to asking the pair
  $J$ and $\Delta$  to form a pseudo-K\"ahler structure.

   \begin{lemma} The pair $(J, \Delta)$ on $\lie h$ forms a pseudo-K\"ahler structure if and only if
   \[
   g(\gamma(x)y,w)-g(\gamma(y)x,w)-g(x,\gamma(y)w)+g(y,\gamma(x)w)=0
   \] for all $x,y,w\in \lie g$.
   \end{lemma}
   \bproof For any $(x,u),(y,v),(z,w)\in \lie g\ltimes V$, expand
   $d\Omega_4((x,u),(y,v),(z,w))$. Since $\gamma$ is torsion-free, it is equal to
   \begin{eqnarray*}
&&-g(\gamma(x)y,w)+g(\gamma(y)x,w)+g(x,\gamma(y)w)-g(y,\gamma(x)w)\\
&&-g(\gamma(z)x,v)+g(\gamma(x)z,v)+g(z,\gamma(x)v)-g(x,\gamma(z)v)\\
&&-g(\gamma(y)z,u)+g(\gamma(z)y,u)+g(y,\gamma(z)u)-g(z,\gamma(y)u)
\end{eqnarray*}
 Since the last three lines are cyclic permutations of $(x,u)$, $(y,v)$ and $(z,w)$, if one of these lines is
 equal to zero, all three equal to zero and therefore $d\Omega_4=0$.  Conversely, if $d\Omega_4=0$,
 set $z=u=v=0$.  Then the last two lines equal to zero, and the lemma follows.
\eproof

Suppose that $(\Omega_c, J)$ is a holomorphic symplectic structure on the semi-direct product
$\lie h=\lie g\ltimes V$ as above. Let $\Omega_3$ and $\Delta$ be the natural
symplectic and pseudo-metric structure on $\lie h$. Assume that $(\Delta, J)$ is pseudo-K\"ahler.
Both $\Omega^{-1}_3$ and $\Omega_4^{-1}$ are candidates to be compatible with $\Lambda=\Omega_c^{-1}$,
so are their linear combinations.
Below is a key technical result in this section.

\begin{proposition}\label{prop:technical}
Suppose that $(\Omega_c, J)$ is a holomorphic symplectic structure on the semi-direct product
$\lie h=\lie g\ltimes V$ as above. Let $\Omega_3$ and $\Delta$ be the natural
symplectic and (pseudo-)metric structure on $\lie h$. Assume that $(\Delta, J)$ is pseudo-K\"ahler structure.
If  there is a real number $\mu$ such that
\begin{equation}\label{constant}
(g^{-1}\omega)(\gamma(a)b)=-4\mu\gamma((g^{-1}\omega)(a))((g^{-1}\omega)(b))
 \end{equation}
 for all $a,b\in \lie g$, then
 \begin{equation}
 \phi=-\frac{i}{4}\Omega_3^{-1}+\mu\Omega_4^{-1}
  \end{equation}
  and $\Lambda=\Omega_c^{-1}$ forms a compatible pair.
\end{proposition}

In the expression (\ref{constant}), we consider the contractions with  $\omega$ and $g$ as
maps from the underlying vector
space $V$ of $\lie g$  to its dual. Therefore, $\jmath=g^{-1}\omega$ is a map from $V$ to $V$.
The following are used frequently in our proof of Proposition \ref{constant} above.
\begin{lemma}\label{lem:technical} Recall that
$\lie h^{1,0}=\{(a,-ia)\in (\lie g\oplus V)_{\CC}: a\in \lie g\}.$
\begin{itemize}
\item As (0,1)-forms, $\Omega_3((a,-ia), -)=-i\Omega_4((\jmath(a), -i\jmath(a)),-).$
\item As (1,0)-forms, $\Omega_3((a,ia), -)=i\Omega_4((\jmath(a), i\jmath(a)),-).$
\item As (1,0)-forms, $\Omega_c((a,-ia), -)=-2\Omega_4((\jmath(a), i\jmath(a)),-).$
\item $\lbra{(a,-ia)}{(b, ib)}^{1,0}=(-\gamma(b)a, i\gamma(b)a)$.
\item $\lbra{(a,-ia)}{(b, ib)}^{0,1}=(\gamma(a)b, i\gamma(a)b)$.
\end{itemize}
\end{lemma}

To prove Proposition \ref{prop:technical}, we consider a generic linear combination of $\Omega_3^{-1}$ and
$\Omega_4^{-1}$, $\phi=\lambda\Omega_3^{-1}+\mu\Omega^{-1}_4$.

Note that we first extend
both $\Omega_3^{-1}$ and
$\Omega_4^{-1}$ by zeroes on $\lie h^{1,0}$. Then they are extended
 as endomorphisms defined on $\lie h^{1,0}\oplus \lie h^{*(0,1)}$
to endomorphisms defined on the exterior
product
$\wedge^\bullet(\lie h^{1,0}\oplus \lie h^{*(0,1)})$ through the identity (\ref{wedge endo}),
by linearity $\phi$ also satisfies (\ref{wedge endo}). Therefore, we will determine the coefficients
$\lambda$ and $\mu$ by solving the non-trivial constraints in (\ref{infinitesimal 1}) and (\ref{phi endo}).

In the current context, the constraint (\ref{infinitesimal 1}) is equivalent to requiring that for all
$\ell_1, \ell_2\in \lie h^{1,0}\oplus \lie h^{*(0,1)}$,
\begin{equation}\label{phi bracket}
\phi(\lbra{\ell_1}{\ell_2})=\lbra{\phi\ell_1}{\ell_2}+\lbra{\ell_1}{\phi\ell_2}.
\end{equation}
Since $\lie h^{1,0}$ is annihilated by $\phi$, and it is closed with respect
to Schouten bracket, if both $\ell_1$ and $\ell_2$ are in $\lie h^{1,0}$, then the identity (\ref{phi bracket})
is trivially satisfied, and hence does not pose any constraint on $\lambda$ and $\mu$.

If $\ell_1\in \lie h^{1,0}$, then there exists $a\in \lie g$ such that $\ell_1=(a,-ia)$.
If $\ell_2\in \lie h^{*(0,1)}$, then there exists $(b, -ib)\in \lie h^{1,0}$ such that
$\ell_2=\Omega_3((b,-ib),-)$. By Lemma \ref{lem:technical},
\[
\ell_2=\Omega_3((b,-ib),-)=\Omega_4((-i\jmath(b), -\jmath(b)), -).
\]
Since $\phi\ell_1=0$,
the constraint in (\ref{phi bracket}) is reduced to
$\phi\lbra{\ell_1}{\ell_2}=\lbra{\ell_1}{\phi\ell_2}$. Since both sides of this identity are
$(1,0)$-vectors, to verify that they are identical, it suffices to show that the evaluation of
any $(1,0)$-forms on these two vectors are identical. Since $\Omega_3$ is non-degenerate, any
$(1,0)$-form has the form
$\Omega_3((n,in), -)$ for some $(0,1)$-vector $(n,in)$. Then a proof of
 (\ref{phi bracket}) is reduced to
check whether the following holds:
\[
\Omega_3((n,in), \phi\lbra{\ell_1}{\ell_2})-\Omega_3((n,in), \lbra{\ell_1}{\phi\ell_2})=0.
\]
Making use of various definitions and Lemma \ref{lem:technical}, we reduce the above identity to
\begin{eqnarray*}
&&\lambda\Omega_3((b,-ib), (\gamma(a)\jmath(n), -\gamma(a)\jmath(n)))
-\lambda\Omega_3((n,in), ([a,b], -i[a,b]))\\
&&+\mu\Omega_4((-i\jmath(b), -\jmath(b)), (i\gamma\jmath(n), -\gamma(a)\jmath(n)))\\
&&-\mu\Omega_4((i\jmath(n), -\jmath(n)), (-i[a, \jmath(b)], -[a, \jmath(b)]))=0.
\end{eqnarray*}
Using definition of $\Omega_3$ and $\Omega_4$ in terms of $\omega$, the above is reduced to
\[
-\lambda\omega(\gamma(b)n, a)+i\mu g(\gamma(\jmath(b))\jmath(n), a)=0.
\]
It is equivalent to
\begin{equation}\label{first constrain}
\lambda \jmath(\gamma(b)n)=i\mu\gamma(\jmath(b))(\jmath(n))
\end{equation}
for all $b, n\in \lie g$. This identity is the first preliminary constraint on $\mu$ and $\lambda$.

Similarly, if $\ell_1, \ell_2\in \lie h^{*(0,1)}$, choose $(a, -ia)$ and $(b,-ib)$ such that
\begin{equation}\label{l 1 2}
\ell_1=\Omega_3((a,-ia), -), \quad \ell_2=\Omega_3((b,-ib), -).
\end{equation}
Since $\lbra{\ell_1}{\ell_2}=0$,
(\ref{phi bracket}) is reduced to
\begin{equation}\label{reduced}
\lbra{\phi\ell_1}{\ell_2}+\lbra{\ell_1}{\phi\ell_2}=0.
\end{equation}
As both terms in the above sum are (0,1)-forms, then its evaluation on any (0,1)-vector $(n, in)$ is equal to
zero. Substitute (\ref{l 1 2}) into identity (\ref{reduced}), evaluate on a (0,1)-vector $(n, in)$,
and make use of Lemma \ref{lem:technical}, we get
\begin{eqnarray*}
&&-\lambda\Omega_3((b,-ib), (\gamma(a)n, i\gamma(a)n))+\lambda\Omega_3((a,-ia), (\gamma(b)n, i\gamma(b)n))\\
&&-\mu\Omega_4((-i\jmath(b),-\jmath(b)), (-i\gamma(\jmath(a))n, \gamma(\jmath(a))n))\\
&&+\mu\Omega_4((-i\jmath(a),-\jmath(a)), (-i\gamma(\jmath(b))n, \gamma(\jmath(b))n))=0.
\end{eqnarray*}

Using definitions of $\Omega_3$ and $\Omega_4$, together with Lemma \ref{lem:technical}, The above
identity is reduced to
\[
\lambda\omega([a,b], n)-i\mu g([\jmath(a), \jmath(b)], n)
\]
for all $n\in \lie g$. That is
\begin{equation}\label{lambda bracket}
\lambda \jmath([a,b])=i\mu [\jmath(a), \jmath(b)].
\end{equation}
Since $\gamma(a)b-\gamma(b)a=[a,b]$ for all $a, b$, the above is equivalent to
\[
\lambda \jmath(\gamma(a)b)-\lambda \jmath(\gamma(b)a)
=-\mu \gamma (\jmath(a))\jmath(b)+\mu\gamma(\jmath(b))\jmath(a).
\]
This identity holds for all $a, b\in \lie g$ so long as
(\ref{first constrain}) holds. Therefore, (\ref{first constrain}) is the only constraint
for solving (\ref{phi bracket}).

Next, we need to find the constraints on $\lambda$ and $\mu$ to satisfy the
identify (\ref{phi endo}). This is equivalent to requiring
\begin{equation}\label{constrain 2}
\lambda(\dbar\Omega_3^{-1}(\ell)-\Omega_3^{-1}\dbar\ell)
+\mu(\dbar\Omega_4^{-1}(\ell)-\Omega_4^{-1}\dbar\ell)=\lbra{\Lambda}{\ell}
\end{equation}
for all $\ell\in \lie h^{1,0}\oplus \lie h^{*(0,1)}$.

Since $\Omega_3^{-1}$ and $\Omega_4^{-1}$ are extended by zero on $\lie h^{1,0}$, when
 $\ell$ is an element in $\lie h^{1,0}$, the constraint (\ref{constrain 2}) is reduced to
 \begin{equation}\label{constrain 2, simplified}
 -\lambda\Omega_3^{-1}\dbar\ell-\mu\Omega_4^{-1}\dbar\ell=\lbra{\Lambda}{\ell}.
 \end{equation}

 Let $A,B$ be elements in $\lie h^{1,0}$, with identity (\ref{21 identity}) and the fact that
  $d\Omega_c=0$, one could check that
\[
\lbra{\Lambda}{\ell}(\Omega_cA, \Omega_cB)=\Omega_c(\ell, \lbra{A}{B}).
\]
If we set $\ell=(x, -ix), A=(a, -ia), B=(b, -ib)$ with $x, a, b\in \lie g$, recall the definitions
 of $\Omega_c$ in terms of $\omega$, then the above is
further simplified to
\begin{equation}\label{lambda ell}
\lbra{\Lambda}{\ell}(\Omega_cA, \Omega_cB)=4i\omega(x, [a,b]).
\end{equation}
In view of (\ref{CR}), the first term on the left-hand-side of the identity in
 (\ref{constrain 2, simplified}) evaluated on
the ordered pair $\Omega_cA, \Omega_cB$ is simplified to
\begin{eqnarray*}
&&-\Omega_3^{-1}\dbar\ell(\Omega_cA, \Omega_cB)\\
&=&-2i\left(\Omega_cB(\lbra{(x, -ix)}{(a,ia)})-\Omega_cA(\lbra{(x,-ix)}{(b, ib)}\right)\\
&=&-2i\left(\Omega_c((b, -ib), \lbra{(x, -ix)}{(a,ia)}^{1,0})
-\Omega_c((a, -ia), \lbra{(x,-ix)}{(b, ib)}^{1,0})\right).
\end{eqnarray*}
With Lemma \ref{lem:technical} and various definitions, one could show that
\begin{equation}\label{Omega 3 without lambda}
-\Omega_3^{-1}\dbar\ell(\Omega_cA, \Omega_cB)=-8\omega(x, [a,b]).
\end{equation}
Similarly,
\begin{eqnarray*}
&&-\Omega_4^{-1}\dbar\ell(\Omega_cA, \Omega_cB)\\
&=&-2\left(
\Omega_cA(-\gamma(\jmath b)x, i\gamma(\jmath b)x)
-\Omega_cB(-\gamma(\jmath a)x, i\gamma(\jmath a)x)
\right).
\end{eqnarray*}
By Lemma \ref{lem:technical}, it is equal to
\[
2\left(
\Omega_4((-2\jmath a, -2i\jmath a), (-\gamma(\jmath b)x, i\gamma(\jmath b)x))
-\Omega_4((-2\jmath b, -2i\jmath b), (-\gamma(\jmath a)x, i\gamma(\jmath a)x))
\right).
\]
By definition of $\Omega_4$, we have
\[
-\Omega_4^{-1}\dbar\ell(\Omega_cA, \Omega_cB)
=8g(i[\jmath a, \jmath b], x).
\]
Since (\ref{first constrain}) is satisfied, (\ref{lambda bracket}) holds. Therefore,
\[
-\mu\Omega_4^{-1}\dbar\ell(\Omega_cA, \Omega_cB)
=8\lambda g(\jmath [a,b], x)=-8\lambda\omega(x, [a,b]).
\]
Combined the above identity with (\ref{Omega 3 without lambda}) and (\ref{lambda ell}),
we obtain
\[
-16\lambda\omega(x, [a, b])=4i\omega(x, [a,b])
\]
for all $x, a, b\in \lie g$.
Therefore, $\lambda=-\frac{i}4$.
Further and similar calculations demonstrate that this is the only constraint \cite{Brian}.

Substitute this constraint into (\ref{first constrain}), we find that $\mu$ is a real number and
for all $a, b\in \lie g$,
\[
\jmath (\gamma(a)b)=-4\mu \gamma(\jmath a)(\jmath b).
\]
It concludes the proof of Proposition \ref{prop:technical}.

\

Let us analyze Proposition \ref{prop:technical} further. If $\mu=0$,  constraint (\ref{constant}) implies
that $\gamma(a)b=0$ for all $a, b\in \lie g$. Therefore, $\gamma=0$. However, the connection $\gamma$ is
torsion-free. This implies that $[a,b]=0$. Therefore, the algebra $\lie h=\lie g\ltimes V$ is trivial.
In particular,  $\Lambda$ is central in the Gerstenhaber algebra $(\wedge^\bullet\lie h, \wedge,
\lbra{-}{-})$, and hence $(\Lambda, \phi=0)$ forms a compatible pair.

Therefore, whenever
$\lie h$ is non-abelian, we may assume that $\mu\neq 0$. In such case, if one multiplies the non-degenerate
bilinear form $g$ on $\lie g$ by the constant $-4\mu$, then the inhomogeneity in equation
(\ref{constant}) allows us to simply this identity to
\begin{equation}
(g^{-1}\omega)(\gamma(a)b)=\gamma((g^{-1}\omega)(a))((g^{-1}\omega)(b)).
 \end{equation}
 Now we could apply  Proposition \ref{prop:technical} and
 Theorem \ref{thm:integrable} to conclude the following.
 \begin{theorem}\label{thm:algebraic}
 Let $\lie g$ be a Lie algebra with an invariant symplectic structure $\omega$ and non-degenerate bilinear
 form $g$. Let $V$ be its underlying vector space. Let $\gamma: \lie g\to \END (V)$ be a torsion-free
 flat connection and $\lie h=\lie g\ltimes_\gamma V$ the associated semi-direct product.
 Then $\lie h$ has a natural complex structure $J$, a symplectic structure $\Omega$
 and a pseudo-metric $\Delta$. If this triple forms a pseudo-K\"ahler structure and if
 \[
(g^{-1}\omega)(\gamma(a)b)=\gamma((g^{-1}\omega)(a))((g^{-1}\omega)(b)).
 \]
 then there exists a deformation from the complex structure $J$ to a symplectic structure $\Omega_2$
 such that $\DGA(J)$ is isomorphic to $\DGA(\Omega_2)$.
 \end{theorem}

\section{Low-dimension examples}\label{sec:example}

According to Andranda \cite{Andranda}, there are three non-trivial four-dimensional complex symplectic algebras.
 Let $e_1, e_2$ be a basis of $\lie g$ and
$v_1, v_2$ be a basis for $V$ such that
\begin{equation}\label{complex in 4d}
Je_1=v_1, \quad Je_2=v_2.
\end{equation}
 Let $e^1, e^2$ and $v^1, v^2$ be the dual bases.
We choose the symplectic structure $\omega$ and the pseudo-metric $g$ on the algebra $\lie g$ to be
\[
\omega=e^1\wedge e^2, \quad g=e^1\otimes e^2+e^2\otimes e^1.
\]
It follows that
\[
\jmath=g^{-1}\omega=e^1\otimes e_1-e^2\otimes e_2.
\]
The natural symplectic form and metric on $\lie g\ltimes V$ are respectively
\[
\Omega=e^1\wedge e^2+v^1\wedge v^2, \quad \Delta=e^1\otimes e^2+e^2\otimes e^1+
v^1\otimes v^2+v^2\otimes v^1.
\]
Moreover, let $z_1=\frac12(e_1-iv_1)$ and $z_2=\frac12(e_2-iv_2)$, $z^1=e^1+iv^1,$ and
$z^2=e^2+iv^2$, then
\begin{eqnarray}
&\Omega_1=-e^1\wedge v^2-v^1\wedge e^2=\frac1{2i}(z^1\wedge z^2-\oz^1\wedge\oz^2)&\label{o 1}\\
&\Omega_2=e^1\wedge e^2-v^1\wedge v^2=\frac12(z^1\wedge z^2+\oz^1\wedge \oz^2)&\label{o 2}\\
&\Omega_3=e^1\wedge e^2+v^1\wedge v^2=\frac12(z^1\wedge \oz^2+\oz^1\wedge z^2)&\label{o 3}\\
&\Omega_4=e^1\wedge v^2-v^1\wedge e^2=\frac{i}2(z^1\wedge \oz^2-\oz^1\wedge z^2)&\label{o 4}
\end{eqnarray}
In particular,
\begin{eqnarray}
\Omega_c=\Omega_1+i\Omega_2=iz^1\wedge z^2,   &\quad&
\Lambda=\Omega_c^{-1}=iz_1\wedge z_2,\nonumber \\
\Omega_3^{-1}=2(z_2\wedge \oz_1+\oz_2\wedge z_1),
 &\quad&
\Omega_4^{-1}=2i(z_2\wedge \oz_1-\oz_2\wedge z_1). \label{o inverse}
\end{eqnarray}

\subsection{Example 1} When the two-dimensional Lie algebra $\lie g$ is
abelian, the only non-trivial object in constructing a four-dimensional semi-direct product in this
case is the torsion-free flat
connection $\gamma$. It is determined by the identities,
\[
\gamma(e_1)v_1=v_2, \quad \gamma(e_1)v_2=0, \quad \gamma(e_2)=0.
\]
Equivalently, the only non-trivial structure equation for $\lie h=\lie g\ltimes V$ is
\[
\lbra{e_1}{v_1}=v_2.
\]
 The dual structure equation is $dv^2=-e^1\wedge v^1.$ Therefore, it is apparent that
$\Omega_4$ is closed, and hence $\lie h$ has a natural pseudo-K\"ahler metric.

As $\jmath e_1=e_1$ and $\jmath e_2=-e_2$,
Proposition \ref{prop:technical} is solved when $\mu=\frac14$.
By the expressions in (\ref{o inverse}),
\[
\phi= -\frac{i}4\Omega_3^{-1}+\frac14\Omega_4^{-1}=\frac14 (\Omega_4^{-1}-i\Omega_3^{-1})=i z_1\wedge \oz_2.
\]
Therefore by Theorem \ref{thm:algebraic},
for the complex structure $J$ in (\ref{complex in 4d})  and the symplectic structure
$\Omega_2$ in (\ref{o 2}), $\DGA(\Omega_2)$ and $\DGA(J)$ are isomorphic and they exist in one
generalized deformation class.

Indeed, for this particular example, the algebraic $\lie h$ is the covering space of the Kodaira-Thurston surface.
It is known that all the concerned cohomology
spaces are given by invariant objects. Therefore, we may also apply
Theorem \ref{thm:weak mirror} on manifold level, and recovers a key result obtained by ad hoc computation in \cite{Poon}.

\subsection{Example 2} In this example, the algebra $\lie g$ is solvable, with structure equation
$[e_1, e_2]=e_2$. The connection $\gamma$ is given by
\[
\gamma(e_1)v_1=-v_1, \quad \gamma(e_1)v_2=v_2, \quad \gamma(e_2)=0.
\]
The structure equations for the semi-direct product $\lie h$ are equivalently given by
\[
de^2=-e^1\wedge e^2, \quad dv^1=e^1\wedge v^2, \quad dv^2=-e^1\wedge v^2.
\]
It follows that $\Omega_4$ is closed.
Further,  $\mu=-\frac14$ solves the constraint in Proposition \ref{prop:technical},
and $\phi=-iz_2\wedge \oz_1$. Therefore, by Theorem \ref{thm:algebraic}
the complex structure $J$ is deformed to a $\Omega_2$ via a holomorphic Poisson structure, and
$\DGA(J)$ is isomorphic to $\DGA(\Omega_2)$.

\subsection{Example 3} In this example, the algebra $\lie g$ is solvable: $[e_1, e_2]=e_2$.
The connection $\gamma$ is given by
\[
\gamma(e_1)v_1=-\frac12 v_1, \quad \gamma(e_1)v_2=\frac12 v_2, \quad \gamma(e_2)v_1=-\frac12 v_2,
\quad \gamma(e_2)v_2=0.
\]
On the semi-direct product the non-trivial structure equations become
\begin{equation}\label{real structure eq}
[e_1, e_2]=e_2, \quad [e_1, v_1]=-\frac12 v_1, \quad
[e_1, v_2]=\frac12 v_2, \quad [e_2, v_1]=-\frac12 v_2.
\end{equation}
The dual equations are
\begin{equation}\label{exterior structure}
de^2=-e^1\wedge e^2, \quad dv^1=\frac12 e^1\wedge v^1, \quad dv^2=-\frac12 e^1\wedge v^2+\frac 12 e^2\wedge v^1.
\end{equation}
It follows that $d\Omega_4=2v^1\wedge e^1\wedge e^2$. In particular, Proposition
\ref{prop:technical} and Theorem \ref{thm:algebraic} are not applicable.
In terms of complex frames, we have
\[
\lbra{z_1}{z_2}=\frac12 z_2, \quad dz^1=-\frac14 z^1\wedge\oz^1,
\quad dz^2=-\frac14(z^1+\oz^1)\wedge z^2-\frac14 z^1\wedge (z^2+\oz^2).
\]
From the differentials, we further obtain that
\begin{equation}
\lbra{z_1}{\oz^1}=\frac14\oz^1, \quad \lbra{z_1}{\oz^2}=-\frac14\oz^2,
\quad \lbra{z_2}{\oz^2}=\frac14\oz^1.
\end{equation}
Taking the complex conjugation, and then the dual expression is
\[
\dbar z_1=-\frac14 \oz^1\wedge z_1-\frac14 \oz^2\wedge z_2,
\quad
\dbar z_2=\frac14 \oz^1\wedge z_2.
\]
As an intermediate step, we put together the structure equation of $\DGA(J)$ on this particular
algebra:
\begin{eqnarray}
&\lbra{z_1}{z_2}=\frac12 z_2, \quad
\lbra{z_1}{\oz^1}=\frac14\oz^1, \quad \lbra{z_1}{\oz^2}=-\frac14\oz^2,
\quad \lbra{z_2}{\oz^2}=\frac14\oz^1&\label{complex structure eq}\\
&\dbar z_1=-\frac14 \oz^1\wedge z_1-\frac14 \oz^2\wedge z_2,
\quad
\dbar z_2=\frac14 \oz^1\wedge z_2,
\quad \dbar\oz^2=-\frac12 \oz^1\wedge \oz^2.&\label{dbar}
\end{eqnarray}

On the other hand,
 \[
 \Omega_2(e_1)=e^2, \quad \Omega_2(e_2)=-e^1, \quad \Omega_2(v_1)=-v^2,
 \quad \Omega_2(v_2)=v^1.
 \]
 Then the linear isomorphism $\Omega_2$ take the Lie bracket on vectors in (\ref{real structure eq}) to
 a Lie bracket on forms.
 \[
 \lbra{e^1}{e^2}=-e^1, \quad \lbra{e^2}{v^2}=\frac12 v^2, \quad
 \lbra{e^2}{v^1}=\frac12 v^1, \quad \lbra{e^1}{v^2}=-\frac12 v^1.
 \]
 With respect to these Lie algebra structures, the first derived subalgebra
 $\lie h^{1,0}\oplus \lie h^{*(0,1)}$  is
 the three-dimensional Heisenberg algebra spanned by $e^1, v^1, v^2$ with $v^1$ being its center.
 In view of the exterior differential as given in (\ref{exterior structure}), $v^1$ is not closed in
  the differential Gerstenhaber algebra of the symplectic structure $\Omega_2$.

  On the other hand, from (\ref{complex structure eq}), we find that the first derived subalgebra
  in $\DGA(J)$ is the three-dimensional Heisenberg algebra spanned by $z_2, \oz^1, \oz^2$ with
  $\oz^1$ being its center. In view of (\ref{dbar}),  $\oz^1$ is $\dbar$-closed.

  Since the center of the derived subalgebra of $\DGA(J)$ is $\dbar$-closed and that of $\DGA(\Omega_2)$ is
  not $d$-closed, then two $\DGA$s could not be quasi-isomorphic

\section*{Remark}

Given the definition of $\Omega_1$ in (\ref{o 1}), it is apparent that $\lie g$ and $V$ are Lagrangian with
respect to $\Omega_1$. As $J\lie g=V$ and $JV=\lie g$, the complex structure $J$ and the complex symplectic
structure is special Lagrangian in the sense of \cite[Definition]{COP}.

Let $\gamma^*$ be the dual representation of $\gamma$, then one obtains the dual semi-direct product
$\widehat{\lie h}=\lie g\ltimes_{\gamma^*} V^*$. Through this space as an intermediate object, it is provided
in \cite[Theorem 5.2]{COP} that there is a natural isomorphism from $\DGA(J)$ to $\DGA(\Omega_1)$. See also
\cite{Ben}.
Therefore, we have
\[
\DGA(\Omega_1)\cong \DGA(J) \cong \DGA(\Omega_2).
\]

\end{document}